\documentclass[table,onecolumn]{IEEEtran}
\usepackage{lipsum,amsmath}
\usepackage{cuted}
\usepackage{bm}

\usepackage{graphicx}
\usepackage{float}
\usepackage{amsmath}
\usepackage{epstopdf}
\usepackage{cases}
\usepackage{color}
\usepackage{algorithm}
\usepackage{booktabs}
\usepackage{algorithmic}
\usepackage{amsmath}
\usepackage{verbatim}
\graphicspath{{figures/}}

\title{Integrated Sensing and Communication enabled Sensing Base Station: System Design, Beamforming, Interference Cancellation and Performance Analysis}
    \author{
        Wangjun Jiang,~\IEEEmembership{Student Member,~IEEE,}
        Zhiqing Wei,~\IEEEmembership{Member,~IEEE,} \\
        Zhiyong Feng,~\IEEEmembership{Senior Member,~IEEE,} \\
        Xu Chen,~\IEEEmembership{Student Member,~IEEE}
        \\
        \thanks{ \emph{Corresponding author: Zhiyong Feng and Zhiqing Wei.}

            W. Jiang, Z. Wei,  Z. Feng, and X. Chen are with the School of Information and Communication Engineering, Beijing University of Posts and Telecommunications, and also with the Key Laboratory of Universal Wireless
            Communications, Ministry of Education, Beijing 100876, China (email: \{jiangwangjun, weizhiqing, fengzy, chenxu96330\}@bupt.edu.cn).

        }
    }

\begin{document}

\maketitle

\begin{abstract}
This paper studies the sensing base station (SBS) that has great potential to improve the safety of vehicles and pedestrians on roads. It can detect the targets on the road with communication signals using the integrated sensing and communication (ISAC) technique.
Compared with vehicle-mounted radar, SBS has a better sensing field due to its higher deployment position, which can help solve the problem of sensing blind areas.
%
In this paper, key technologies of SBS are studied, including the beamforming algorithm, beam scanning scheme, and interference cancellation algorithm.
To transmit and receive ISAC signals simultaneously, a double-coupling antenna array is applied.
The free detection beam and directional communication beam are proposed for joint communication and sensing to meet the requirements of beamwidth and pointing directions.
The joint time-space-frequency domain division multiple access algorithm is proposed to cancel the interference of SBS, including multiuser interference and duplex interference between sensing and communication.
Finally, the sensing and communication performance of SBS under the industrial scientific medical power limitation is analyzed and simulated.
Simulation results show that the communication rate of SBS can reach over 100 Mbps and
the range of sensing and communication can reach about 500 m.
\end{abstract}


\section{Introduction}\label{Introduction}

\subsection{Background and Motivation}
The causes of traffic accidents have shown that 90-93\% of vehicle incidents are caused by human errors \cite{[Driverless_Future]}. If there were no human error, accidents on the road would have been reduced by more than 90\% \cite{[Driverless_Future]}. The cooperative intelligent transportation system (C-ITS) is one of the solutions to reduce human error \cite{[C-ITS]}. A complete C-ITS requires both sensing and communication functions \cite{[RSU]}. {\color{black} Sensing in C-ITS} is mainly realized by sensors, such as radar installed on vehicles or roadside units (RSUs) \cite{[RSU]}.
Compared with vehicle-mounted sensing equipment (VMSE) \cite{[VMSE]}, RSUs have a better sensing field due to their higher deployment position, which can help solve the problem of sensing blind areas \cite{[RSU_VMSE]}.

In current C-ITS, roadside sensing and communication devices are individually designed \cite{[C-ITS]}.
In order to fully obtain the advantages of the integrated sensing and communication (ISAC) technique, this paper puts forward the scheme of the sensing base station (SBS).
Different from the traditional RSU, SBS realizes radar sensing by exploiting the communication signal using the ISAC technique, which allows sensing and communication systems to share spectrum bands, and can thus improve spectrum utilization.
Although there have been relevant research on SBSs, and few papers have conducted detailed research on specific key technologies of SBSs.
Recognize this fact, key technologies of SBS are studied in this paper, including the beamforming algorithm  \cite{[JWJ]}, beam scanning scheme, and interference cancellation algorithm, which will be introduced in Section \ref{sec:system-model} and \ref{sec:avoiding-interference}.

\subsection{Related Work} \label{sec:Introduction-2}

RSU is a roadside equipment that can support high-speed vehicle to infrastructure (V2I) communication and sensing services \cite{[RSU_Energy]}. Liu {\it{et al.}} proposed a kind of RSU that can sense road condition information through a 79 GHz millimeter radar and other sensing equipment \cite{[RSU]}.
Compared with the above studies, SBS proposed in this paper realizes the sensing function with the mobile communication signal, which can improve the efficiency of spectrum and energy resources, and reduce the system size and the system cost \cite{[JRC_1]}.

As wireless devices and data traffic grow exponentially, spectrum is becoming increasingly scarce \cite{[spectrum_limit]}. {\color{black} Based on ISAC,} the sensing and communication subsystems can realize spectrum sharing between sensing and communication. {\color {black} The ISAC technology} has great application potential in vehicular networks \cite{[JRC_AV]}, where communication signals can also be used for environmental sensing. In the field of the 6th generation (6G) wireless communication system, ISAC is an increasingly important technology to communicate and sense the environment \cite{[6G]}.

With the similarity of radio frequency front-end architecture between wireless communication and radar sensing, ISAC technology has been developed rapidly \cite{[JRC_waveform]}. The carrier frequency used in communication systems has shifted to the microwave regime (300 MHz $ \sim $ 3 THz) that is traditionally used for radar. Therefore, ISAC is a feasible solution to realize spectrum sharing between sensing and communication, so as to alleviate the scarcity of spectrum resources.
There are some related studies in ISAC signal design, such as orthogonal frequency division multiplexing (OFDM) signal \cite{[OFDM]}, linear frequency modulation (LFM) \cite{[LFM]} and phase-modulated continuous wave (PMCW) \cite{[PMCW]}.
In terms of communication, OFDM signal has significant advantages of overcoming multipath interference and exploiting frequency diversity \cite{[OFDM_2]}. It is the signal model of the 4th generation (4G) and the 5th generation (5G) mobile communication system.
In terms of radar sensing, OFDM signal has significant advantages of high precision and high resolution \cite{[OFDM_advantage]}. The number and frequency interval of sub-carriers can be adjusted flexibly to obtain the ambiguity function of ``pushpin shape" \cite{[OFDM_3]}.
Therefore, we choose the OFDM signal as the signal model of SBS.

SBS with ISAC technology faces a great challenge of interference cancellation. Compared with the traditional base station, the interference cancellation of SBS is more complex because SBS needs to deal with not only multiuser interference but also duplex interference between sensing and communication.
For multiuser interference, Sit {\it{et al.}} proposed a signal separation algorithm based on continuous wave OFDM signal reconstruction \cite{[IIE-1],[IIE-2],[IIE-3]}. However, the performance of this algorithm is greatly affected by carrier frequency offset estimation. We adopt a multiple access scheme using space, time and frequency resources to cancel multiuser interference.
Duplex interference between sensing and communication is caused by SBS's echo signal and received uplink communication signal.
There are many duplex modes in mobile communication, such as time division duplex (TDD) and frequency division duplex (FDD). The TDD mode has significant advantages of high spectral efficiency, and it is the main duplex mode of the current 5G system. Therefore, we present a scheme based on TDD to cancel the duplex interference between sensing and communication.

\subsection{Contributions of Our Work}
In this paper, we focus on studying SBS that has great potential to improve the safety for vehicles and pedestrians on roads. Contributions of this paper are summarized as follows.

1. We propose the system model of SBS that can detect the targets on road with communication signals using ISAC technique.
To realize the ISAC technique, we adopt a double-coupling antenna array to transmit and receive ISAC signals simultaneously.
Moreover, key technologies of SBS are studied, including antenna array design, beamforming algorithm, beam scanning scheme, and interference cancellation algorithm.

2. The free detection beam (FDB) and directional communication beam (DCB) are proposed for joint communication and sensing to meet the requirements of beamwidth and pointing directions. Although J. Andrews {\it {et al.}} proposed the beamforming algorithm for FDB and DCB in \cite{[Multibeam]}, which is to optimize FDB and DCB respectively, and then superpose the weight vectors of FDB and DCB linearly. When the number of beams increases, sidelobe leakage will be superposed, thus reducing the gain of the main lobe. For point-to-point interaction scenarios, only one DCB is needed, and the problem of sidelobe leakage is not severe. For SBS, multiple communication beams are needed to provide communication services for multiple users simultaneously. To solve this problem, we propose an improved beamforming algorithm to jointly optimize the weight vectors of FDB and DCBs. Simulation results verify the advantages of the proposed beamforming algorithm. Details are introduced in section \ref{sec:BF}.

3. We analyze the interference of a single SBS system, including multiuser interference and duplex interference between sensing and communication. A signal frame structure is designed to cancel duplex interference.
The joint time-space-frequency domain division multiple access (TSF-DMA) algorithm is proposed to cancel the multiuser interference. Details are introduced in section \ref{sec:avoiding-interference}.

4. The sensing and communication performance of SBS under the industrial scientific medical power limitation is analyzed and simulated. We derive the expression of the accuracy of ranging and velocity estimation, and carry out the simulation verification. Details are introduced in Section \ref{sec:Performance} and \ref{sec:Numerical results}.

The remaining parts of this paper are organized as follows.
Section \ref{sec:system-model} describes the SBS system, followed by the design of antenna arrays, ISAC signal model, beamforming algorithm and beam scanning scheme. Section \ref{sec:avoiding-interference} introduces a scheme to cancel the interference of SBS, including multiuser interference and duplex interference between sensing and communication. The interference cancellation algorithm is also introduced in this section. Section \ref{sec:Performance} provides a detailed derivation of the communication and sensing performance of SBS. In section \ref{sec:Numerical results}, we present simulation results of the SBS system performance. Section \ref{sec:Conclusion} concludes the paper.

The symbols used in this paper are described as follows. Vectors and matrices are denoted by boldface small and capital letters; the transpose, complex conjugate, Hermitian, inverse, and pseudo-inverse of the matrix ${\bf{A}}$ are denoted by ${{\bf{A}}^{T}}$, ${{\bf{A}}^*}$, ${{\bf{A}}^{H}}$, ${{\bf{A}}^{ - 1}}$ and ${{\bf{A}}^\dag}$, respectively.

\section{SBS System}\label{sec:system-model}

In this section, we will introduce the system model of SBS, followed by key technologies of SBS, including the design of antenna arrays, ISAC signal model, beamforming algorithm and beam scanning scheme.

\subsection{System Model}

\begin{figure}[ht]
	\includegraphics[scale=0.65]{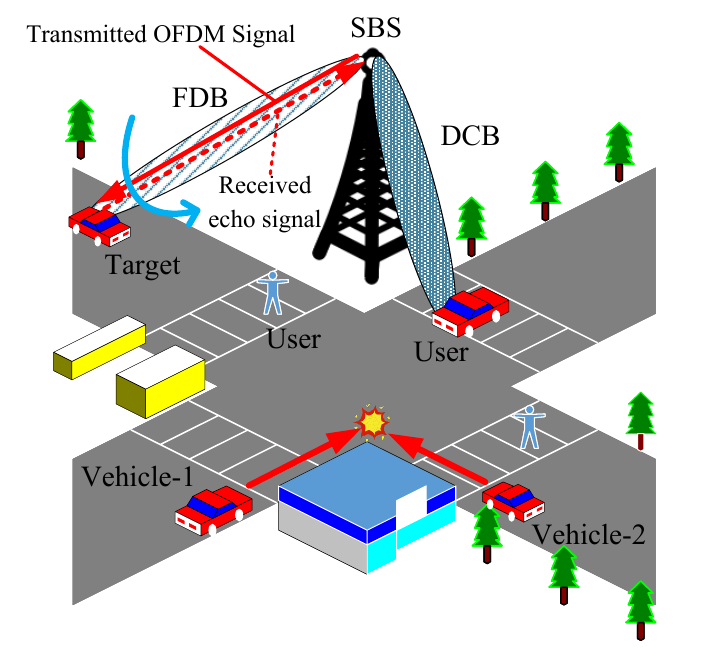}
	\centering
	\caption{The SBS forms FDB and DCB to detect targets and communicate with users.}
	\label{fig:SBS_system}
\end{figure}

As Fig. \ref{fig:SBS_system} shows, the SBS is located on the roadside.
Because of the building blockage, vehicle-1 and vehicle-2 {\color{black} cannot} detect each other by VMSE. As a result,
vehicle-1 and vehicle-2 have a risk of collision.
The SBS installed at the top of the road can detect the two vehicles in advance and transmit the sensing information to the two vehicles by ISAC signals, further helping them avoid collision.
{\color {black}
The detected targets can be divided into two parts: cooperative targets, and non-cooperative targets. Cooperative targets are users that can communicate with SBS, such as vehicles. Non-cooperative targets cannot communicate with SBS, such as trees.
For non-cooperative targets, the SBS can send the location information of non-cooperative targets to users to help users be informed of their environment and thus avoid collisions.
For cooperative targets ({\color{black} users}), the SBS needs to associate the user's location information with the digital identity.
In this paper, we focus on the non-cooperative targets' detection to help users be informed of their environment and thus avoid collisions.
}

{\color{black}
Sensing can be divided into two stages: target search, and target tracking. Target search requires SBS to generate an omnidirectional scanning beam to detect the target, and target tracking requires SBS to generate a directional beam directed to the target.
Communication can also be divided into two stages: initial user access, and user directed communication. User initial access requires SBS to generate an omnidirectional scanning beam to try to establish a communication connection with the user, and user directed communication requires SBS to generate a directional beam directed to the user.
Thus, both the sensing service for target search and the communication service for initial user access can be realized by the ISAC signals carried by the FDB. Both the sensing service for target tracking and the communication service for user directed communication can be realized by the ISAC signals carried by the DCB.
In this paper, we focus on the sensing service for target search and the communication service for user directed communication.
}

\begin{figure}[ht]
	\includegraphics[scale=0.4]{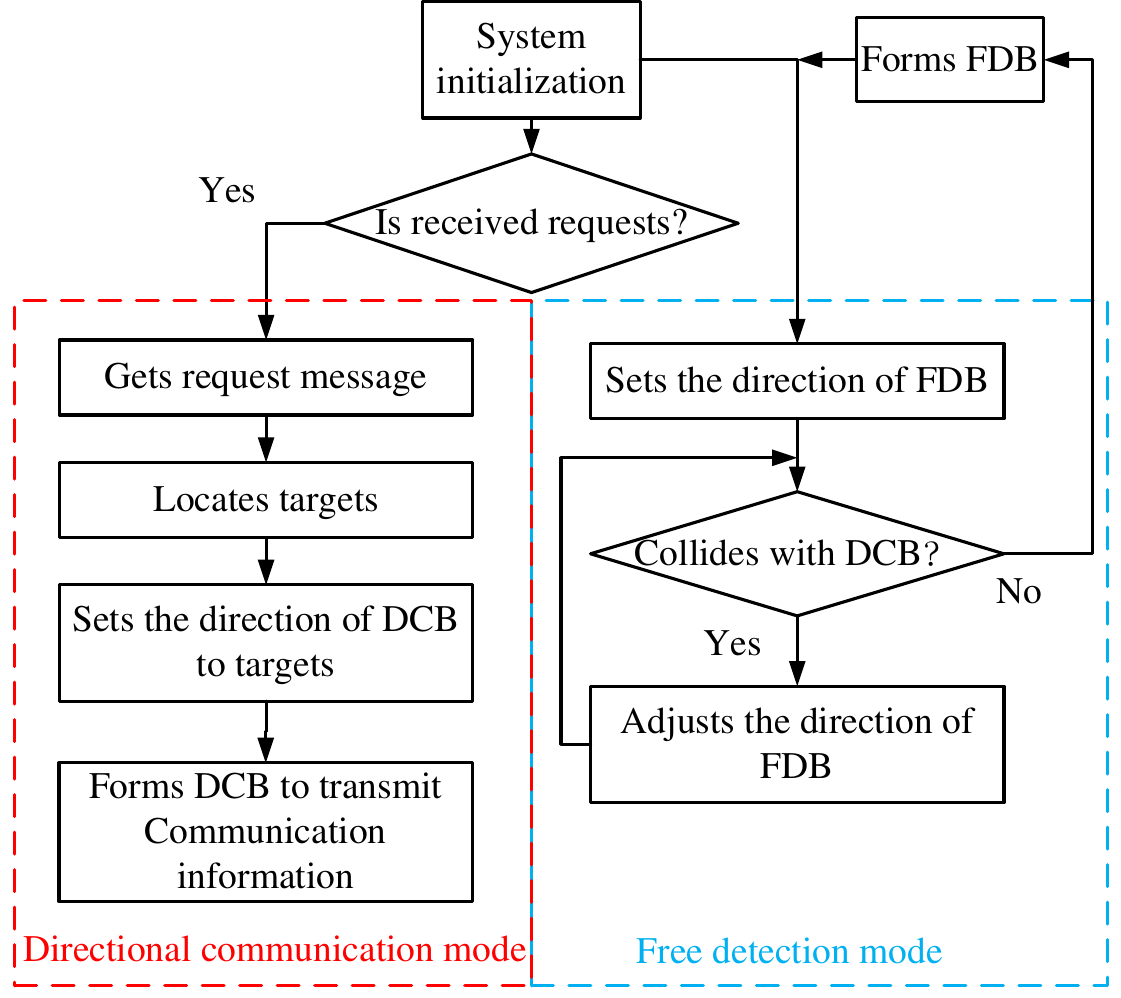}
	\centering
	\caption{Process of the SBS system.}
	\label{fig:SBS_process}
\end{figure}

To provide communication and sensing services for vehicles, the SBS system has two working modes: free detection mode and directional communication mode.
As Fig. \ref{fig:SBS_process} shows, SBS will enter free detection mode after system initialization. If SBS receives communication requests from users, it will enter directional communication mode.

$\cdot$ In free detection mode, SBS will form a FDB to detect targets by ISAC signals. On the one hand, SBS can get the location of targets by the ISAC echo signals.
{\color{black}
On the other hand, SBS can provide the omnidirectional scanning communication service to users by the ISAC signals.
}

$\cdot$ In directional communication mode, SBS has located users in free detection mode, it can form a DCB with high beam gain pointing at users, which can improve the signal-to-noise ratio (SNR) of received signals, thus improving the communication performance.

$\cdot$ It is noted that the free detection mode and directional communication mode are running in parallel.

%

\subsection{OFDM Signal Model} \label{subsec: OFDM_Signal}

Since OFDM signals have advantages in sensing and communication introduced in section \ref{sec:Introduction-2}, the OFDM signal is chosen as the signal model of SBS, which can be expressed as \cite{[OFDM]}
{\color{black}
\begin{equation}\label{equ:OFDM_signal}
	\begin{aligned}
	x(t) = \sum\limits_{m = 0}^{M - 1} \sum\limits_{n = 0}^{N - 1} & {s\left( mN + n \right)} \cdot \\
	& {\rm exp} \left({j2\pi{f_n}t} \right) \cdot {\rm rect} \left (\frac{{t - mT}}{T} \right),
	\end{aligned}
\end{equation} }
where $M$ is the number of OFDM symbols, $N$ is the number of subcarriers, $m$ is the OFDM symbol index, $n$ is the index of subcarriers, $s\left(mN + n\right)$ is the complex modulation symbol, ${\rm rect} \left(\cdot \right)$ is the rectangle function. The bandwidth of subcarriers is ${f_n} = {n \cdot \Delta f}$, where ${\Delta f}$ denotes subcarrier spacing. The duration time of each OFDM symbol $T$ contains the elementary symbol duration ${T_{os}}$ and the guard interval ${T_g}$. If ${T_g}$ is larger than the maximum multipath delay, inter-symbol interference will be eliminated. In the reception and processing of OFDM sensor signals, the signal-based sensor processing method can use the modulation symbol domain to estimate the doppler frequency shift and delay of received signals. Details of this method are described in section \ref{sec:sensing-Performance-2}.

\subsection{Channel Model} \label{sec:communication-Performance-1}
The Rayleigh fading channel model and Rician fading channel model are both classical wireless communication channel models. The Rayleigh fading channel model is used to describe multipath channels and is suitable for scenarios without line of sight (LOS) paths. On the contrary, the Rician fading channel model is suitable for scenarios with a LOS path, such as vehicular communications \cite{[rice_vehicular]}. Communication between SBS and users is realized by DCB and exists a LOS path. Hence, the Rician fading channel model is adopted as the channel model of the SBS system. Then, the power of the communication signal received by SBS is
\begin{equation}\label{equ:rician}
	\begin{split}
		{P_{r,c}} = \frac{\rho{P_{t}} {G_{{\rm{com}}}} {\it{g}_{p,c}} \xi \lambda^2}{\left(4\pi\right)^2{x}^{\alpha }}
	\end{split},
\end{equation}
where ${P_t}$ denotes the transmitting signal power, $\lambda$ is the wavelength of the communication signal, $g_{p,c}$ is the communication processing gain \cite{[OFDM]}, ${G_{{\rm{com}}}}{\rm{ = }}{g_{t}} \cdot  {g_{r,c}}$ is the product of transmitting beam gain and receiving beam gain, $\alpha$ is the large scale path loss exponent, and $x$ is the distance between SBS and users.
{\color{black}
The power ratio for communication, $\rho $ indicates the distribution of the communication and sensing power. When $\rho = 0$, the SBS is fully in sensing mode and all the transmit power is used for sensing. When $\rho = 1$, the SBS is fully in communication mode, and all transmitting signal power is used for communication, which is the same as the conventional base station communication system.}
{\color{black}Assuming that the perfect channel state information (CSI) is available at the receiver, the probability density function (PDF) of the received SNR, denoted by $\xi $, can be written as  \cite{[RICE],[Xinyuan],[Xinyuan-4],[Xinyuan-4-31]}  }
\begin{equation}\label{equ:rician_distribution}
	\begin{split}
		{f_\xi }\left( w \right) = \frac{{\left( {K + 1} \right)}}{{\bar \xi }{e^{{  \frac{{\left( {K + 1} \right)w}}{{\bar \xi }}} +  K}}}{I_0}\left( {2\sqrt {\frac{{K\left( {K + 1} \right)w}}{{\bar \xi }}} } \right)
	\end{split},
\end{equation}
where ${I_0}\left(\cdot \right)$ is the zero-order modified
Bessel function of the first kind. The statistically average SNR $\bar \xi $ can be expressed as
\begin{equation}\label{equ:average_SNR}
	\begin{split}
		\bar \xi \rm{=} \frac{\rho{\it{P}_{t}} {\it{G}_{{\rm{com}}}} \it{g}_{p,c} \lambda^2} {\left(4\pi \right)^2{\it{x}}^{ \alpha}\it{P}_n\it{F}_n}
	\end{split},
\end{equation}
where $\it{P}_n$ and $\it{F}_n$ denote the noise power and noise figure, respectively.
The Rician factor ${\it{K}}=\frac{{{{\it{v}}^2}}}{{{\varsigma ^2}}}$ is defined to be the ratio between the signal power of the LOS path, denoted by ${v^2}$, and the multipath reflected
power is denoted by ${\varsigma ^2}$.

\subsection{Design of SBS Antenna Array}

\begin{figure}[ht]
	\includegraphics[scale=0.4]{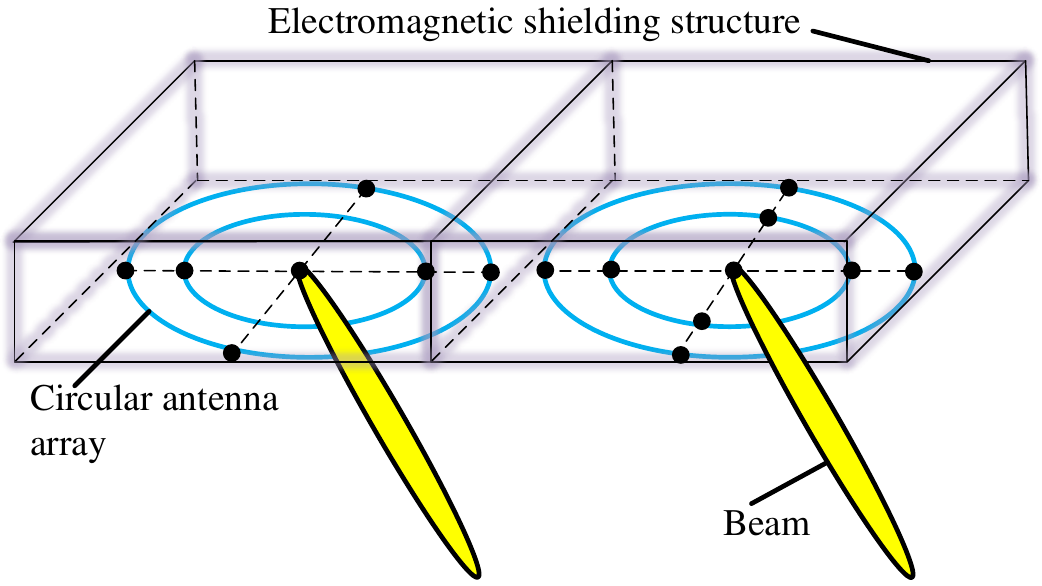}
	\centering
	\caption{Two decoupled circular uniform subarrays.}
	\label{fig:SBS_ant_1}
\end{figure}

\textbf{}
\begin{figure}[ht]
	\includegraphics[scale=0.7]{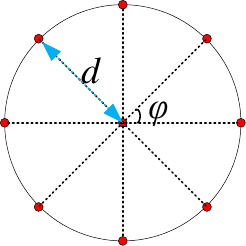}
	\centering
	\caption{Single layer of SBS antenna array.}
	\label{fig:SBS_ant_3}
\end{figure}

Since targets and users are distributed randomly on the road, the antenna array installed at SBS needs to form three-dimensional (3D) downward beams. The circular uniform antenna array is a feasible choice to generate 3D beams \cite{[ant_design]}.
{\color {black}
Compared with traditional uniform planar arrays, circular uniform antenna arrays can generate a beam with smaller sidelobe \cite{[UCA]}.
}
As Fig. \ref{fig:SBS_ant_1} shows, the SBS antenna (SBSA) consists of two decoupled circular uniform subarrays. The electromagnetic shielding structure is adopted to weaken the coupling interference between SBSAs. There are ${\it{p}}$ circles in each SBSA. We define each circle as a layer.
Regarding the center antenna as the phase reference antenna (PFA), {\color{black} the index of layers from center to periphery is from ${{0}}$ to ${{\it{p}}-1}$. Except for the ${{0}}$-th layer with only one antenna element, there are ${\rm{2}}^b$ antenna elements in each layer, where ${\it{b}}$ is an integer}. Antenna elements in each layer are located on the evenly distributed polar angle, as shown in Fig. \ref{fig:SBS_ant_3}. The distance between the antenna elements located at the same angle of adjacent layers is ${\it{d}}$.

To avoid phase ambiguity, the distance between adjacent antenna elements should meet the following conditions
\begin{equation}\label{equ:distance_of_antenna}
	\begin{aligned}
		{\it{d}} &\le {\frac{\lambda}{2}} \\
		2{\it{d}}\sin \frac{\varphi }{2} &\le \frac{\lambda}{2}
	\end{aligned},
\end{equation}
where $\varphi = \frac{2 \pi}{2^b}$ is the angle difference of adjacent elements within the same layer. It indicates that both the distance between adjacent layers and the distance between adjacent elements within the same layer should be smaller than half of the wavelength.

\section{Key technologies of SBS} \label{sec:SBS-key}

\subsection{Beamforming Algorithm for FDB and DCB} \label{sec:BF}

\begin{figure}[ht]
	\includegraphics[scale=0.55]{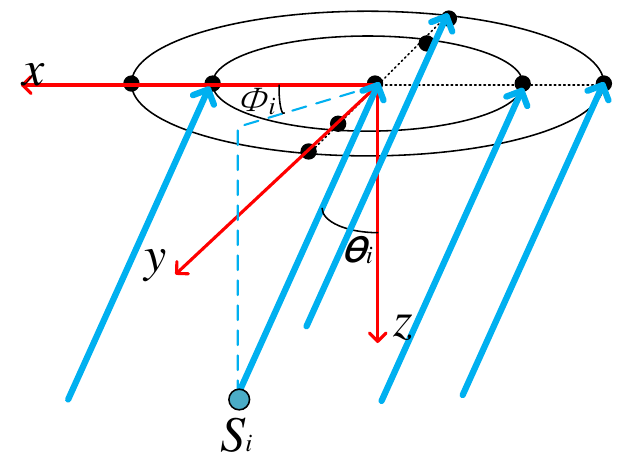}
	\centering
	\caption{Receiving signals of SBSA.}
	\label{fig:SBS_ant_2}
\end{figure}

As Fig. \ref{fig:SBS_ant_2} shows, assuming that there are ${\it{l}}$ far field signals arriving at SBSA from different directions. The $\it{i}$-th planar far field signal comes from the direction of $\left({\phi _i},{\theta _i}\right)$. The phase difference between the ${\it{n}}$-th antenna element of the ${\it{m}}$-th layer and PFA caused by the ${\it{i}}$-th planar far field signal can be expressed as
\begin{equation}\label{equ:daoxiang}
	{{\it{a}}_{m,n}}\left({\phi _i},{\theta _i}\right) = {\rm exp} \left ( - j\frac{{2\pi }}{\lambda } {\bf{q}} _{m,n}^T{ {\bf{v}} _i} \right ), {\color{black} \in {{\mathcal C}^{1 \times 1}}}
\end{equation}
where the polar mapping vector ${{\bf{v}} _i}$ can be described as
\begin{equation}\label{equ:polar_mapping_vector}
	{{\bf{v}} _i} = {\left[{\rm cos} {\phi _i}{\rm sin}{\theta _i},{\rm sin}{\phi _i}{\rm sin}{\theta _i}\right]^T} . {\color{black} \in {{\mathcal C}^{2 \times 1}}}
\end{equation}
The distance between PFA and the ${\it{n}}$-th antenna element of the ${\it{m}}$-th layer can be expressed as
\begin{equation}\label{equ:distance_between_array}
	{\bf{q}} _{m,n} = {\left[{\rm cos} \left({\psi _{m,n}} \right)md,{\rm sin} \left({\psi _{m,n}} \right)md \right]^T}, {\color{black} \in {{\mathcal C}^{2 \times 1}}}
\end{equation}
where ${\psi _{m,n}}{{ = n}} \cdot \varphi $ denotes the polar angle of ${\it{n}}$-th antenna element of ${\it{m}}$-th antenna layer.
The steering vector of signal ${{\it{s}}_{\it{i}}}$ with angle of arrival $\left({\phi _{\it{i}}},{\theta _{\it{i}}}\right)$ can be expressed as
{\color{black}
\begin{equation}\label{equ:steering_vector}
	\begin{aligned}
	 {\bf{ a}}_{\it i} & =
	  [  \rm{1},{\it{a}_{\rm{1,0}}}\left({\phi _{\it i}},{\theta _{\it i}}\right), \cdots, \\
	  & {\it{a}_{{\rm{1},{\rm{2}^{\it b}} - \rm{1}}}}\left({\phi _{\it i}},{\theta _{\it i}}\right),
		\cdots,{\it{a}_{\it{p} - \rm 1,{\rm{2}^{\it b}} - 1}}\left({\phi _{\it i}},{\theta _{\it i}}\right) ]  ^{\it T}, {\color{black} \in {{\mathcal C}^{N_t \times 1}}}
	\end{aligned}
\end{equation} }
where $\color{black} N_t = \left(p-1\right)2^b + 1$ is the number of antenna elements, $\rm{1}$ denotes that the phase of PFA is zero and the remaining items in \eqref{equ:steering_vector} refer to the phase difference between other antenna elements and PFA.
{\color{black}
Then, the expected direction set of generated beam can be expressed as $\mathcal{S}_a = \left\{  \left({\theta _{\rm 1}}, {\phi _{\rm 1}}\right) , \cdots, \left({\theta _{\it i}}, {\phi _{\it i}}\right) , \cdots, \left({\theta _{\it l}}, {\phi _{\it l}}\right)           \right\}$
}
If there are ${\it{l}}$ far field signals, then the steering matrix ${\bf{D}}$ can be described as
\begin{equation}\label{equ:steering_matrix}
	{\bf{D}} = \left[ {{{\bf a}_1},{{\bf a}_2},...{{\bf a}_{\it i}}...,{{\bf a}_l}} \right], {\color{black} \in {{\mathcal C}^{N_t \times l}}}
\end{equation}
where ${{\bf a}_{\it i}}$ is the steering vector of the ${\it i}$-th signal.
{\color{black}
After beamforming, the received signal vector for SBSA can be expressed as
\begin{equation}\label{equ:signal_after_beamforming}
	\begin{aligned}
		{y}  &= {{\bf{w}}^H} \left({\bf{D}} {{\bf x}}  + {{\bf n}}\right), {\color{black} \in {{\mathcal C}^{1 \times 1}}}
	\end{aligned},
\end{equation}
where ${{\bf x} = \left[x_1, x_2, \cdots , x_l \right]^T} {\color{black} \in {{\mathcal C}^{l \times 1}}} $ is the signal vector, $\bf {{n}} $ is the additive gaussian white noise vector.
The weight vector of each antenna, ${\bf{w}}$ can be expressed as
\begin{equation}\label{equ:weight_vector}
	{{\bf{w}}} = \left[ {{{w}}_{0,0}}{\rm{,}}{{{w}}_{1,0}},{{{w}}_{1,1}},...{{{w}}_{{\it {m,n}}}}...,{{{w}}_{\it p - \rm 1,{2^{\it b} - \rm 1}}} \right]^T, {\color{black} \in {{\mathcal C}^{N_t \times 1}}}
\end{equation}
where ${{{w}}_{{\it {m,n}}}}$ is the weight of the $\it n$-th antenna element of the $\it m$-th layer.
}

The main challenge of beamforming for ISAC signals is that communication and sensing have different requirements of beamwidth and pointing directions \cite{[Multibeam]}. Sensing requires fast direction switching FDBs to obtain fast environment sensing, while communication requires stable and accurately-pointed DCBs to obtain high-capacity communication.
{\color{black}
To meet the different requirements of sensing and communication, the beamforming algorithm for FDB and DCB is proposed in \cite{[Multibeam]}. Two beamforming algorithms are proposed in \cite{[Multibeam]}: local optimization and global optimization. The core idea of local optimization is to optimize FDB and DCB respectively, and then superpose the weight vectors of FDB and DCB linearly. The core idea of global optimization is to find optimal solutions for a single beamforming vector.
Moreover, \cite{[Multibeam]} only considers the scenario with a FDB and a DCB, the sidelobe leakage of FDB  and DCB is not serious. However, multiple communication beams are needed for SBS to provide communication services for multiple users simultaneously. When the number of beams increases, sidelobe leakage will be superposed.
%
To solve this problem, we jointly optimize the weight vectors of FDB and DCBs.
\subsubsection{Problem formulation} \label{sec:BF-1}
In this paper, we intend to minimize the error between the desired amplitude directional response and the generated amplitude directional response
\begin{equation}\label{equ:convex-1}
\min \limits_{{{ {\bf{w}} }_{opt}}}  \begin{Vmatrix}
	 { {\bf{w}}^H _{{{opt}}}}{\bf D}_d - {{\bf r}_{ad}}
\end{Vmatrix}^2,
\end{equation}
where
\begin{itemize}
	\item ${\bf D}_d {\color{black} \in {{\mathcal C}^{N_t \times N_d}}} $ is the steering matrix of all possible directions, {\color{black} with $N_d \ge l$ being the number of all possible  directions of the beam,}
	\item ${\bf{w}}_{opt} {\color{black} \in {{\mathcal C}^{N_t \times 1}}} $ is the final beamforming vector,
	\item ${{\bf r}_{ad}} {\color{black} \in {{\mathcal C}^{1 \times N_d}}}$ is the desired amplitude directional response
	\begin{equation}\label{equ:r_ad}
		{{\bf{ r}}_{ad}}\left(\phi ,\theta \right)=\left\{
		\begin{aligned}
			1, & \quad \left(\phi,\theta \right) {\rm{ \in  \mathcal{S}_a}} \\
			0, & \quad \left(\phi,\theta \right) {\rm{ \notin  \mathcal{S}_a}}
		\end{aligned},
		\right.
	\end{equation}
which means that the generated beam has an amplitude response of 1 in the desired direction and 0 in the other directions.
\end{itemize}
If the problem \eqref{equ:convex-1} is solved directly without constraint, ${\bf{w}}_{opt}$ may grow excessively, resulting in distortion of the generated beam and large sidelobe.
Regularization is a typical method to solve this problem, and $L_2$ regularization is adopted in this paper. Then the optimization problem can be described as
\begin{equation}\label{equ:convex-2}
	\min \limits_{{{ {\bf{w}} }_{opt}}}  \begin{Vmatrix}
		{ {\bf{w}}^H _{{{opt}}}}{\bf D}_d - {{\bf r}_{ad}}
	\end{Vmatrix}^2 + \beta \begin{Vmatrix}
	{\bf{w}}_{opt}
\end{Vmatrix}^2,
\end{equation}
where
\begin{itemize}
	\item $\beta $ is the regularization factor, which is set as 0.01.
\end{itemize}
However, it is difficult to obtain two types of FDB and DCB by directly solving \eqref{equ:convex-2}, so we use the beamforming algorithm in \cite{[Multibeam]} to obtain the beamforming vector of single beam. Then, the final beamforming vector can be expressed as
\begin{equation}\label{equ:convex-3}
{ {\bf{w}} _{{opt}}} =   {\bf W} {{\bf f} _w}, {\color{black} \in {{\mathcal C}^{N_t \times 1}}}
\end{equation}
where
\begin{itemize}
	\item $\bf{W}$ is the weight matrix of FDB and DCBs
	\begin{equation}\label{equ:W}
			{\bf{W}} = \left[ {\bf w}_{\it f}{\rm{,}}{{\bf {w}}_{{{1}}}}{\rm{,}}...{\rm{,}}{{\bf {w}}_{{d}}} \right], {\color{black} \in {{\mathcal C}^{N_t \times \left(d+1\right)}}}
	\end{equation}
	\item  ${\bf w}_{\it f}$ and ${\bf w}_{i, i=1,...,d}$ are the weight vectors of FDB and DCBs, which can be obtained by the beamforming algorithm in \cite{[Multibeam]},
	\item $d$ is the number of DCBs,
	\item ${ {\bf{f}} _w} = \left[ {{\it{f}}_f}{\rm{,}}{{\it{f}}_1}{\rm{,}}{{\it{f}}_2}{\rm{,}}...{\rm{,}}{{\it{f}}_d}\right]^T {\color{black} \in {{\mathcal C}^{\left(d+1\right) \times 1}}}$ is a vector of combining coefficients for $\bf{W}$.
\end{itemize}
So far, the optimization problem of beamforming can be expressed as
\begin{equation}\label{equ:convex}
	\begin{aligned}
		 \mathop {\min }\limits_{{{ {\bf{f}} }_w}} &
		\begin{Vmatrix}
			{ {\bf{w}}^H _{{{opt}}}}{\bf D}_d - {{\bf r}_{ad}}
		\end{Vmatrix}^2
		+ \beta
		\begin{Vmatrix}
			{{\bf w} _{{{opt}}}}
		\end{Vmatrix}^2,\\
		& {\rm{s.t: }} { {\bf{w}} _{{opt}}} = {\bf W} {{\bf f} _w}.
	\end{aligned}
\end{equation}
%
%

{\color{black}
\subsubsection{Problem Solution} \label{sec:BF-2}
The problem \eqref{equ:convex} is a quadratic optimization problem under linear constraints, which is a convex optimization problem. We can solve the problem by using MATLAB’s cvx toolbox.
The pseudo-code for solving \eqref{equ:convex} is shown in Algorithm \ref{alg:CVX}.
\begin{algorithm}[htb]
	\caption{Pseudo-code for solving \eqref{equ:convex} by cvx toolbox}
	\label{alg:CVX}
	\begin{algorithmic}
		\STATE cvx\_begin
		  \STATE \qquad variable ${\bf f}_w$
		  \STATE \qquad constant ${\bf W}$
		  \STATE \qquad constant ${\bf r}_{ad}$
		  \STATE \qquad constant $\beta$
		  \STATE \qquad minimize $\begin{Vmatrix}
		  	{ {\bf{w}}^H _{{{opt}}}}{\bf D}_d - {{\bf r}_{ad}}
		  \end{Vmatrix}^2
		  + \beta
		  \begin{Vmatrix}
		  	{{\bf w} _{{{opt}}}}
		  \end{Vmatrix}^2$
		  \STATE \qquad \quad subject to ${ {\bf{w}} _{{opt}}} =  {\bf W} {{\bf f} _w}$
		 \STATE cvx\_end
	\end{algorithmic}
\end{algorithm}
}

Details of the beamforming {\color{black} algorithm} is shown in Algorithm \ref{alg:BF}. }
	\begin{algorithm}[htb]
		\caption{Beamforming algorithm of FDB and DCBs}
		\label{alg:BF}
		\begin{algorithmic}
			\STATE $\bf Input$:
			The regularization factor $\beta $.
			\STATE $\bf Output$:
			Final beamforming vector: ${\bf{w}} _{{opt}}$. \\
			\STATE $\bf 1)$	Generate the desired amplitude directional response, ${\bf r}_{ad}$ based on \eqref{equ:r_ad};
			\STATE $\bf 2)$ Obtain the weight vectors of FDB, ${\bf w}_{f}$, and DCBs, ${\bf w}_{i}$, based on the beamforming algorithm proposed in \cite{[Multibeam]}, and construct the weight matrix, $\bf W$ based on \eqref{equ:W};
			\STATE $\bf 3)$ Construct the convex optimization problem \eqref{equ:convex} of combining the weight vectors of FDB, ${\bf w}_{f}$, and DCBs, ${\bf w}_{i}$;
			\STATE $\bf 4)$ Solve the convex problem by CVX toolbox to obtain the combining coefficients, ${\bf f}_{w}$;
			\STATE $\bf 5)$ Reconstruct the final beamforming vector of FDB and DCBs based on \eqref{equ:convex}.
		\end{algorithmic}
	\end{algorithm}
\subsection{Scanning Mode of FDB}

\begin{figure}[ht]
	\includegraphics[scale=0.75]{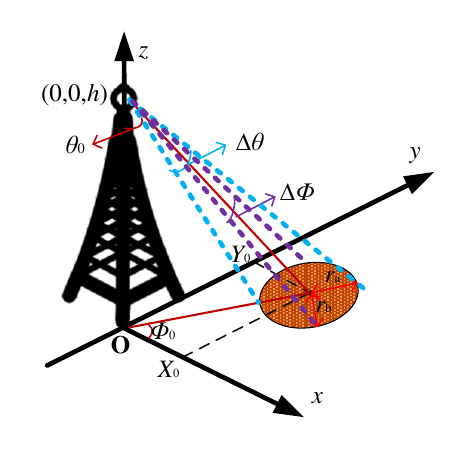}
	\centering
	\caption{Coverage of FDB on the road.}
	\label{fig:SBS_beam_coverage}
\end{figure}

\begin{figure}[ht]
	\includegraphics[scale=0.65]{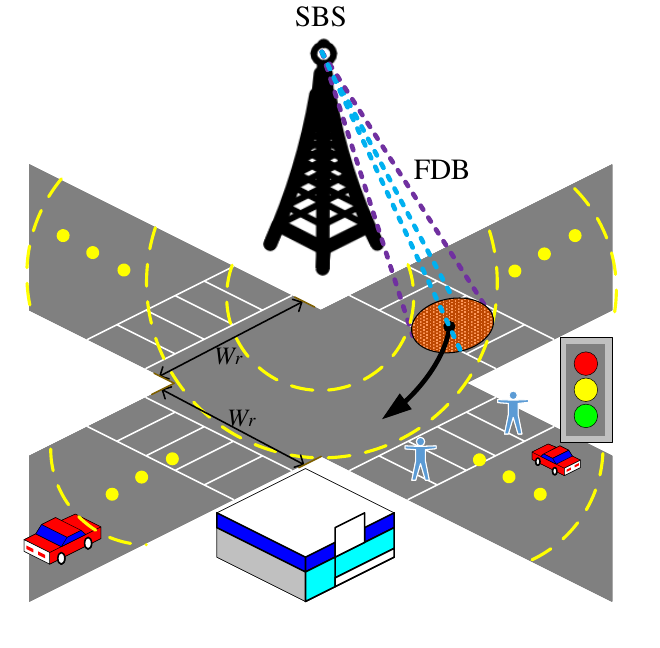}
	\centering
	\caption{Scanning mode of FDB.}
	\label{fig:SBS_scan_mode}
\end{figure}

As illustrated in Fig. \ref{fig:SBS_beam_coverage}, the sensing unit generated on the cross-traffic road region by conical FDB is elliptic, and the location of the sensing unit $(X_0,Y_0)$ can be derived as
\begin{equation}\label{equ:beam_location}
	\left(X_0,Y_0 \right) = \left({h \cdot {{{\tan}}  \left( {\theta _0} \right)  } {{\cos}} \left( {\phi _0} \right)  }  ,{h \cdot {{{\tan}}  \left(  {\theta _0} \right) } {{\sin}} \left( {\phi _0} \right) }  \right),
\end{equation}\\
where $\it {h}$ is the height of SBS, $\left({\phi _0},{\theta _0}\right)$ is the direction of FDB.
The semi-major axis ${{\it{r}}_a}$ and semi-minor axis ${{\it{r}}_b}$ of elliptical sensing unit can be expressed as
\begin{equation}\label{equ:beam_coverage}
	\begin{aligned}
		{r_a} & = \frac{{ {\it{h}} \cdot \tan \left ( \Delta \theta \right ) }}{{\cos \left ({\theta _0} \right )}} \cdot \frac{1}{{\cos \left({\theta _0}\right)}} \approx \frac{{{\it{h}} \cdot \Delta \theta }}{{{{\cos }^2} \left ({\theta _0} \right )}},\\
		{r_b} &= \frac{{h \cdot {{\tan}} \left(\Delta \phi \right)}}{2{{\rm{cos}} \left({\theta _0} \right)}},
	\end{aligned}
\end{equation}\\
where $\Delta \theta $ and $\Delta \phi $ are the pitch width and the azimuth width of FDB, respectively.

In order to sense targets on the road, FDB needs to scan in the cross-traffic road region.
As Fig. \ref{fig:SBS_scan_mode} shows, the cross-traffic road region can be expressed as
\begin{equation}\label{equ:road_region}
	\begin{aligned}
		A_r = \{ \left(x,y\right),  & \left[ \left(0 \le x \le W_r \right) \cup \left(-W_r \le y \le 0\right) \right] \\
		& \cap  \left(x^2+y^2+h^2 \le {\color {black} R^2_{\max,s} }\right)  \},
	\end{aligned}
\end{equation}
where $W_r$ is the width of the road, $R_{\max,s}$ is the maximum sensing range that is derived in section \ref{sec:sensing-Performance-1}.
In order to avoid interference with the signals carried by FDB and DCB, FDB should skip the region of DCB.
The complete scanning mode is introduced in algorithm \ref{alg:Framwork_3}, where $\tau$ is the beam dwell duration, and $T_{sc}$ is the scanning period of FDB.
\begin{algorithm}[htb]
	\caption{Scanning period of FDB}
	\label{alg:Framwork_3}
	\begin{algorithmic}
		\STATE
		$\bf Input$: $\it{h}$, $W_r$, ${R _{\max,s }}$ and the beamwidth, $\left(\Delta \phi ,\Delta \theta \right)$.
		\STATE
		$\bf Output$: The scanning period, ${T_{sc}}$.
		\STATE $\bf 1)$ Initial the starting direction of FDB, $\left(\phi_0 ,\theta_0 \right) = \left(0,\Delta \theta \right)$. \\
		\STATE $\bf 2)$ Initial the location of the sensing unit, $\left(X_0 ,Y_0 \right)$, which is obtained by \eqref{equ:beam_location}. \\
		\STATE $\bf 3)$ Initial the beam dwell duration, $\tau = 10$ ms.
		\STATE $\bf 4)$ Initial the scanning period of FDB, ${T_{sc}}=0$.
		\WHILE{${X_0}^2+{Y_0}^2+{\it{h}}^2 \le {R^2 _{\max,s }}$}
		\STATE Step 1: Update $\left(\phi_0 ,\theta_0 \right)$ as follows.
		\STATE $\phi_0  = \phi_0  + \Delta \phi $.
		\IF {$\phi_0  \ge 360^\circ$}
		\STATE $\theta_0  = \theta_0  + \Delta \theta $.
		\STATE $\phi_0  = \phi_0  - 360^\circ$.
		\ENDIF
		\STATE Update $\left(X_0,Y_0\right)$ according to \eqref{equ:beam_location}.
		\IF {$\left(X_0,Y_0\right) \in A_r$ }
		\STATE Step 2: Avoid collision between  FDB and DCB.
		\IF {DCB is not in the direction of $\left(\phi_0 ,\theta_0 \right)$}
		\STATE SBS Forms FDB in the direction of $\left(\phi_0 ,\theta_0 \right)$.
		${T_{sc}}={T_{sc}}+\tau $.
		\ENDIF
		\ENDIF
		\ENDWHILE
	\end{algorithmic}
\end{algorithm}
%
%

\subsection{Interference Cancellation of SBS}\label{sec:avoiding-interference}

In this section, we will introduce the interference of single SBS system, which consists of multiuser interference and duplex interference between sensing and communication.

\subsubsection{Duplex Interference between Sensing and Communication}\label{sec:avoiding-interference-1}
{
\color{black}
Duplex interference consists of two main components: interference of the sensing echo signal to the uplink communication signal and interference of the transmitted sensing signal to the uplink communication signal.

$\cdot$ Interference of the sensing echo signal to the uplink communication signal: when the receiving antenna receives both the sensing echo signal and the uplink communication signal, the two will interfere with each other. To address this problem, the signal frame structure is proposed to stagger the sensing echo signal and the uplink communication signal by time division.

$\cdot$ Interference of the transmitted sensing signal to the uplink communication signal: Due to the use of continuous waves for sensing, SBS needs to control two antenna arrays at the same time to achieve simultaneous transmission of the sensing signal and reception of the echo signal.
Due to the close proximity of the transmitting antenna array and the receiving antenna array, the transmitting signal will leak to the receiving antenna array, which will cause the interference between sensing and communication. Physical isolation of antenna arrays and self-interference signal processing are generally used to mitigate this type of interference. Related work can be referred to \cite{[duplex_interference]}.
}

\begin{figure}[ht]
	\includegraphics[scale=0.6]{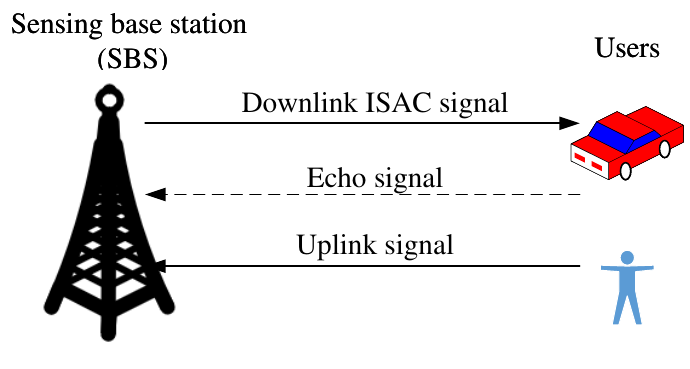}
	\centering
	\caption{Duplex interference between sensing and communication.}
	\label{fig:SBS_interference}
\end{figure}

\begin{figure}[ht]
	\includegraphics[scale=0.4]{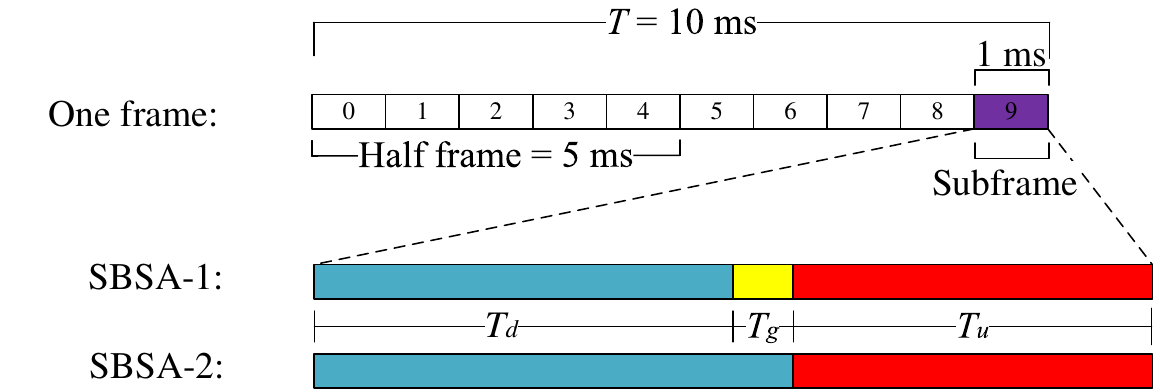}
	\centering
	\caption{Signal frame structure.}
	\label{fig:SBS_time_resource}
\end{figure}
As Fig. \ref{fig:SBS_interference} shows, echo signals of downlink ISAC signals may interfere with the uplink signal. We propose a signal frame structure to cancel the duplex interference.
As Fig. \ref{fig:SBS_time_resource} shows, each frame consists of two half frames with a period of 5 ms. Each half frame is divided into three parts downlink interval (DI) for ${{\it{T}}_d}$ ms, guard interval (GI) for ${{\it{T}}_g}$ ms and uplink interval (UI) for ${{\it{T}}_u}$ ms.

The DI is designed for sensing and downlink communication.
In DI, SBSA-1 and SBSA-2 are in different states. SBSA-1 is sending the downlink ISAC signal, while SBSA-2 is receiving the ISAC echo signal.
The GI is designed to avoid interference between the uplink communication signal and echo signal.
In GI, SBSA-1 is in  the idle state, while SBSA-2 can keep receiving the ISAC echo signal.
The UI is designed for uplink communication.
In UI, SBSA-1 and SBSA-2 are both in the receiving state for receiving the uplink signal.
Detail configuration of DI, GI and UI is listed in table \ref{label:time_resource}.

\begin{table}[ht]
	\caption{Configuration of DI, GI and UI}
	\label{label:time_resource}
	\begin{tabular}{l|l|l|l}
		\hline \hline
		& DI & GI & UI                                                      \\ \hline
		SBSA-1 & \begin{tabular}[c]{@{}l@{}}Transmits \\ downlink \\ ISAC signal\end{tabular} & In idle state           & \begin{tabular}[c]{@{}l@{}}Receives \\ uplink signal\end{tabular} \\ \hline
		SBSA-2 & \multicolumn{2}{l|}{Receives ISAC echo signal}   & \begin{tabular}[c]{@{}l@{}}Receives \\ uplink signal\end{tabular}   \\ \hline
	\end{tabular}
\end{table}

\subsubsection{Multiuser Interference}\label{sec:avoiding-interference-2}
The multiple access scheme is a feasible scheme to cancel multiuser interference. Orthogonal frequency division multiple access (OFDMA) is a main multiple access methods of current mobile communications. OFDMA has significant advantages of high spectrum efficiency and anti-multipath fading \cite{[OFDMA]}. Since SBS communicates
with users by forming DCB, space resources can be used for multiuser interference cancellation. We adopt the TSF-DMA algorithm to cancel multiuser interference.

\begin{figure}[ht]
	\flushleft
	\includegraphics[scale=0.6]{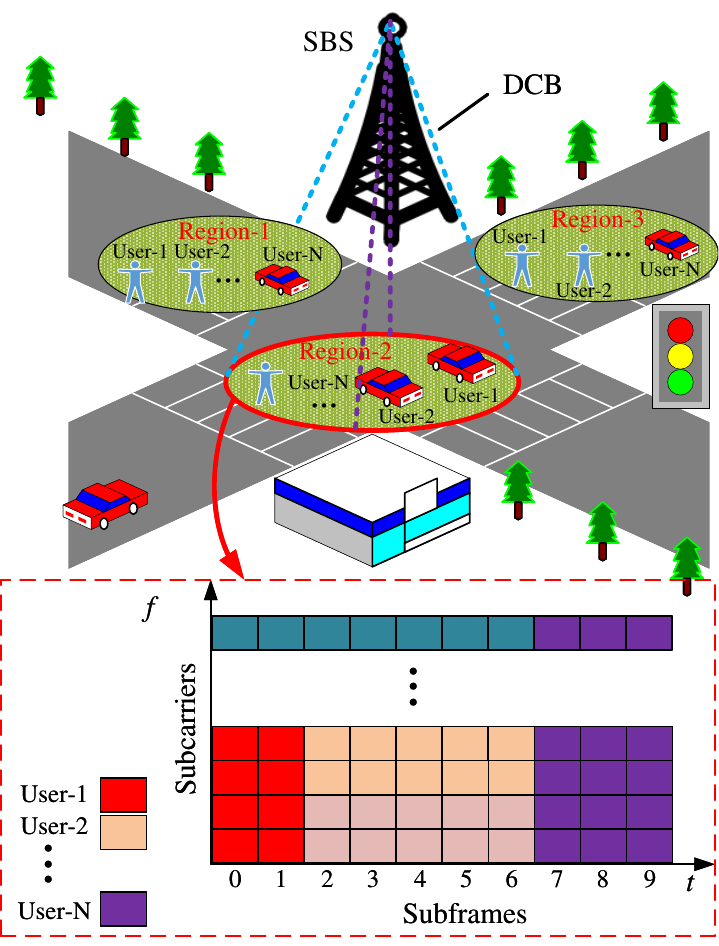}
	\caption{The TSF-DMA scheme of SBS.}
	\label{fig:SBS_resource}
\end{figure}

As Fig. \ref{fig:SBS_resource} shows, on the one hand, different users in different DCBs can be distinguished by different space resources. On the other hand, different users in the same DCB can be distinguished by different time and frequency resources.
In the SBS system, the smallest time domain resource is subframe, which is to satisfy the signal frame structure mentioned above.

\section{Performance Analysis of SBS}\label{sec:Performance}

\subsection{Communication Performance of SBS}\label{sec:communication-Performance}

\subsubsection{Outage Probability} \label{sec:communication-Performance-2}

Outage probability is the probability that the instantaneous SNR at the receiver is below a threshold $\xi_{th}$ required for successful reception \cite{[Outage_Pro]}, which is
\begin{equation} \label{equ:interrupt}
	\begin{split}
		P_c^{out} = {\Pr} \left( {\frac{{{P_{r,c}}}}{\it{P}_n\it{F}_n} < \xi_{th} } \right)
	\end{split},
\end{equation}
where $P_n$ is the noise power, $\xi_{th}$ is acceptable SNR threshold, which depends on quality of service (QoS) requirements,  ${\Pr}\left(\cdot\right)$ denotes probability. Substitute \eqref{equ:rician} into \eqref{equ:interrupt}, the outage probability can be derived as
\begin{equation}\label{equ:interrupt_2}
	\begin{split}
		P_c^{out} = \Pr \left( {\xi  < \frac{{ \left(4\pi\right)^2\xi_{th} \it{P}_n\it{F}_n{x^\alpha }}}{{\rho{{{P}}_{t}}{G_{{\rm{com}}}}\it{g}_{p,c}{\lambda^{\rm{2}}}}}} \right)
	\end{split}.
\end{equation}
Substitute \eqref{equ:rician_distribution} and  \eqref{equ:average_SNR} into \eqref{equ:interrupt_2}, $P_c^{out}$ can be derived as \cite{azari2018ultra}
\begin{equation}
	\begin{aligned}
		& P_c^{out}\left( {x,\xi_{\it{th}} } \right)
		\rm{=} 1 - \\
		& \qquad Q\left( {\sqrt {2K} ,\sqrt {\frac{{2 \cdot \left(4\pi\right)^2\xi_{th} \left( {1 + K} \right){\it{x}^\alpha }\it{P}_n\it{F}_n}}{{\rho{{{\it{P}}}_{t}}{\it{G}_{{\rm{com}}}}\it{g}_{p,c}{\lambda^{\rm{2}}}}}} } \right),
	\end{aligned}
	\label{equ:interrupt_3}
\end{equation}
where $Q\left( {{q_1},{q_2}} \right)$  is the first order Marcum Q-function. Noting that the relationship between outage probability $P_c^{out}$ and successful transmission probability $P_c^s$ is
\begin{equation}
	\begin{aligned}
		P_c^s &= 1 - P_c^{out} \\
		&= Q\left( {\sqrt {2K} ,\sqrt {\frac{{2 \cdot \left(4\pi\right)^2\xi_{th} \left( {1 + K} \right){\it{x}^\alpha }\it{P}_n\it{F}_n}}{{\rho{{{\it{P}}}_{t}}{\it{G}_{{\rm{com}}}}\it{g}_{p,c}{\lambda^{\rm{2}}}}}} } \right).
	\end{aligned}
	\label{equ:ProbSuccBorder}
\end{equation}

\subsubsection{Maximum Communication Range} \label{sec:communication-Performance-3}

In order to ensure the quality of communication, the outage probability should be smaller than threshold $\varepsilon$ \cite{[Outage_Pro]}. The maximum communication range $R_{\max,c}$ is defined as the distance between SBS and users when the outage probability $P_c^{out} = \varepsilon$. When $\it{x} < R_{\max,c}$, the outage probability $P_c^{out} < \varepsilon$. According to \eqref{equ:interrupt_3}, the maximum communication range is
{
\color{black}
\begin{equation}
	\begin{aligned}
	& R_{\max,c} \rm{=} \\
	& \left (  \left( {{\it{Q}}^{{\rm{ - }}1}}(\sqrt{2\it{K}},1-\varepsilon)\right)^2 \cdot  \frac{\rho\it{P}_t\it{G}_{\rm{com}}\it{g}_{p,c}{\lambda^{\rm{2}}}}{2 \cdot (4\pi)^2\xi_{th}(1+\it{K})\it{P}_n\it{F}_n} \right )^\frac{1}{\alpha},
	\end{aligned}
	\label{equ:Communication_range}
\end{equation}
}
where ${Q^{{\rm{ - }}1}}$ is the inverse function of Marcum Q-function. Since the first order Marcum Q-function is a monotone decreasing function of the second parameter ${p_2}$ \cite{[Q_func]}, the inverse function of Marcum Q-function denoted by ${Q^{{\rm{ - 1}}}}\left( {q_1,q_2} \right)$ can make the expression $Q\left( {q_1,{Q^{{\rm{ - 1}}}}\left( {q_1,q_2} \right)} \right) = q_2$ valid.

\subsubsection{Outage Capacity} \label{sec:communication-Performance-4}

The outage capacity ${C_{out}}$ can be treated as the throughput of the typical H-O pair link which is the communication performance of SBS \cite{[Outage_Pro]}.
The outage capacity ${C_{out}}$ can be expressed as
\begin{equation}
	{C_{out}} {\rm{ = }} B \cdot P_c^{s} \cdot \log_2 \left( {1 + {\xi_{th} }} \right),
	\label{equ:OutThrou}
\end{equation}
where $ B $ is the total signal bandwidth.
Based on \eqref{equ:ProbSuccBorder}, and \eqref{equ:OutThrou}, the outage capacity at different distances can be expressed as
{\color{black}
\begin{equation}
	\begin{aligned}
		{C_{out}}\left( x \right) & {\rm{ = }} {B} \cdot \log_2 \left( {1 + {\xi_{th} }} \right) \cdot \\
		& Q\left( {\sqrt {2K} ,\sqrt {\frac{{2 \cdot (4\pi)^2\xi_{th} \left( {1 + K} \right){\it{x}^\alpha }\it{P}_n\it{F}_n}}{{\rho{{{\it{P}}}_{t}}{\it{G}_{{\rm{com}}}}\it{g}_{p,c}{\lambda^{\rm{2}}}}}} } \right).
	\end{aligned}
	\label{equ:OutThrouResult}
\end{equation}
}

\subsection{Sensing Performance of SBS}\label{sec:sensing-Performance}

{ The sensing area size of SBS is limited by the maximum sensing range, and the accuracy of range and velocity estimation is also an important indicator of SBS sensing ability. In this section, the maximum sensing range of SBS and the accuracy of range and velocity estimation are analyzed respectively. At the same time, the ability of SBS to quickly sense the environment is limited by the scanning period of FDB, which will also be analyzed in this section.}

\subsubsection{Maximum Sensing Range} \label{sec:sensing-Performance-1}
According to \cite{[OFDM]}, the maximum sensing range is
{\color{black}
\begin{equation}
	R_{\max,s} = \left (\frac{(1-\rho) {P_{t}} {{\rm{g}}_{t}} {{\rm{g}}_{r,s}} {{\rm{g}}_{p,s}} \sigma{\lambda^2} }{ {{{{(4\pi )}^3} \cdot ({P_n}{F_n}+{P_I}) \cdot {\it \Gamma}_{\min}} }} \right )^{\frac{1}{4}},
	\label{equ:Radar_range}
\end{equation}
}
where $\sigma $ is the sensing cross section, ${{\it{g}}_{r,s}}$ is the sensing receiving beam gain, $P_I$ is the ground clutter power, ${\it \Gamma}_{\min}$ is the minimum required signal-to-interference-plus-noise ratio (SINR) for efficient sensing. The sensing processing gain ${{\it{g}}_{p,s}=M \cdot N}$ is introduced in \cite{[OFDM]}.

\subsubsection{Accuracy of ranging and velocity estimation} \label{sec:sensing-Performance-2}
	Since the range and velocity obtained through the OFDM integrated signal are discrete and have quantization errors \cite{[OFDM]}, most studies adopt the resolution to measure the detection capability, and few studies on the accuracy are made. Li {\it{et al}}. has derived a rough representation of ranging and velocity accuracy based on OFDM integrated signal, and proposed two algorithms to improve the accuracy: fractional Fourier transformation (FRFT) and phase analysis \cite{[OFDM_resolution]}.
	{\color{black} Considering the concrete white gaussian noise and compound-gaussian distribution ground clutter model, the accuracy of ranging and velocity measurement is deduced in this section.}

	Based on the OFDM signal model mentioned in section \ref{subsec: OFDM_Signal}, the received echo signal matrix ${\bf{S}}_{R}$ can be expressed as
	\begin{equation}\label{equ:OFDM_receive_all}
		{\bf{S}}_{R} = \begin{bmatrix}
			{s_{Rx}}(0,0)& \cdots  & {s_{Rx}}(0,N-1) \\
			{s_{Rx}}(1,0)& \cdots  & {s_{Rx}}(1,N-1) \\
			\vdots & \ddots & \vdots \\
			{s_{Rx}}(M-1,0)& \cdots  & {s_{Rx}}(M-1,N-1)
		\end{bmatrix}.
	\end{equation}
	The received modulation symbol ${s_{Rx}}(m,n)$ can be expressed as \cite{[OFDM]}
	\begin{equation}\label{equ:OFDM_receive}
		\begin{aligned}
			{s_{Rx}}(m,n) &= A_s\left(m,n\right) {s_{Tx}}\left(m,n\right) {\exp} \left( - j2\pi {{\rm{f}}_n}\frac{{2{R_r}}}{c}\right) \\
			&  \cdot {\exp} \left( - \frac{j4\pi {\it{m}}T{{V_r}{f_c}}}{c} \right) + A_n + A_I ,
		\end{aligned}
	\end{equation}
	where the complex amplitude factor $A_s(m,n)$ denotes the attenuation and phase shift occurring due to the propagation and scattering process, $P_s(m,n) = \begin{Vmatrix}
		A_s(m,n)
	\end{Vmatrix}^2$ is the power of each received modulation symbol, ${s_{Tx}}(m,n)$ is the sending modulation symbol, $c$ is the speed of light, $f_c$ is carrier frequency, ${R_r}$ and $V_r$ are the range and velocity of the detected target relative to SBS respectively. The rectangular window function with time delay can be neglected in baseband processing \cite{[OFDM]}.
	The additive white compound-gaussian noise $A_n = a_n + ib_n$, where $a_n$ and $b_n$ follow gaussian distribution $ {{\mathcal{CN}}\rm{(0,}}\sigma _n^2{\rm{)}}$, the probability density function is
	\begin{equation}
		{\it{f}_{Ni}(a_n)} = \frac{1}{{\sqrt{2\pi}{\sigma_n}}}\exp{\left( - \frac{a_n^2}{{2\sigma_n^2}}\right)}.
		\label{equ:noise}
	\end{equation}
	The amplitude and phase of noise $A_n$ follow Rayleigh distribution and uniform distribution respectively. The power of noise is $P_n = 2\sigma^2_n$.
	The compound-gaussian ground clutter $A_I = a_i + ib_i$, where $a_i$ and $b_i$ follow gaussian distribution $ {{\mathcal{CN}}\rm{(0,}}\sigma _i^2{\rm{)}}$,
	The amplitude and phase of the ground clutter $A_I$ follow Rayleigh distribution and uniform distribution respectively \cite{[rayleigh]} \cite{[compound-gaussian]}. The power of the ground clutter is $P_I = 2\sigma^2_i$.
	Let us define $N_I = A_n + A_I$, which follows compound-gaussian distribution. The amplitude and phase of $N_I$ follow Rayleigh distribution and uniform distribution respectively.

	After removing the transmitted information from the received information symbols by an element-wise complex division, the elements of division matrix expression can be derived as
	\begin{equation}\label{equ:OFDM_signal_1}
		\begin{aligned}
			{\left({{\bf{S}}_g}\right)_{m,n}} &= \frac{\left({{\bf{S}}_{Rx}}\right)_{m,n}}{\left({{\bf{S}}_{Tx}}\right)_{m,n}} \\
			&=
			{A_s}(m,n){k_r(n)}{k_v(m)} + \frac{N_I}{{s_{Tx}(m,n)}}
		\end{aligned},
	\end{equation}
	where
	\begin{equation}\label{equ:OFDM_signal_4}
		{\left({{\bf{S}}_{Tx}}\right)_{m,n}} = {{{{\rm{s}}_{Tx}}(mN + n)}},
	\end{equation}
	\begin{equation}\label{equ:OFDM_signal_5}
		{\left({{\bf{S}}_{Rx}}\right)_{m,n}} = {{{{\rm{s}}_{Rx}}(mN + n)}},
	\end{equation}
	\begin{equation}\label{equ:OFDM_signal_2}
		{k_r(n)} = {\rm exp} \left(\frac{ - j4\pi {n \Delta f}{{R_r}}}{c} \right),
	\end{equation}
	\begin{equation}\label{equ:OFDM_signal_3}
		{k_v(m)} = {\rm exp} \left( - \frac{j4\pi {\it{m}}T{{V_r}{f_c}}}{c} \right),
	\end{equation}
	Assuming that OFDM signals adopt the quadrature amplitude modulation (QAM) mode, then the amplitude of sending modulation symbols is the same. Assuming that each received modulation symbol has the same power. Then ${{\bf{S}}_g}$ can be simplified as
	\begin{equation}\label{equ:OFDM_signal_new_4}
		\begin{aligned}
			{\left({{\bf{S}}_g}\right)_{m,n}} &= {A_s}{k_r(n)}{k_v(m)} + N_I.
		\end{aligned}
	\end{equation}
	Applying doppler Fourier transformation (DFT) for each row of ${{\bf{S}}_g}$, the velocity of target ${V_r}$ can be deduced as follows
	\begin{equation}\label{equ:Doppler}
		{V_{r}} \in \left [\frac{c \cdot {in{d_{s,n}}}}{{2TMf_c}},\frac{c \cdot ({in{d_{s,n}} + 1})}{{2TMf_c}} \right),
	\end{equation}
	where $in{d_{s,n}}$ is the index of the peak of the DFT of the ${\it{n}}$-th row of ${{\bf{S}}_g}$, can be expressed as
	\begin{equation}\label{equ:Doppler_ind_1}
		{in{d_{s,n}}} = {\max}  \left\{v(0,n),v(1,n), \cdots ,v(M-1,n) \right\},
	\end{equation}
	where
	\begin{equation}\label{equ:Doppler_ind_2}
		\begin{aligned}
			v(i,n) &= \mathsf {DFT} \left[{A_s}{k_r(n)}{k_v(m)}+{A_n e^{j \varphi_n} + A_I e^{j \varphi_I}}\right] \\
			&= {A_s}{k_r(n)} \sum_{m=0}^{M-1} k_v(m)e^{-j2\pi\frac{im}{M}} + N_Ie^{-j2\pi\frac{im}{M}}
		\end{aligned}.
	\end{equation}
	Based on \eqref{equ:Doppler}, we can get the velocity of target estimated by this algorithm \cite{[OFDM]} is discrete and has quantization error. Averaging multiple values is an effective algorithm to eliminate quantization error. Therefore, without loss of generality, we can neglect the quantization error and select the real velocity of the target as $V_r = \frac{k \cdot c}{2f_cTM}$, which $k$ is an integer, then $k_v(m) = e^{-j2\pi \frac{mk}{M}}$.
	Then, the probability of estimated velocity of the target $V_e = \frac{k \cdot c}{2f_cTM}$ can be derived as
	\begin{equation}\label{equ:Doppler_ind_3}
		\begin{aligned}
			p\left(V_e = \frac{k \cdot c}{2f_cTM}\right) &= \prod_{i=0,i\ne k}^{M-1} p\left(v(k,n) >  v(i,n)\right) \\
			&= \left(1-e^{-\frac{2\pi^2 A_s^2 k_r^2(m)}{\sigma^2_n + \sigma^2_i}}\right)^{(M-1)} \\
			&= \left(1-e^{-{4\pi^2\gamma k_r^2(m) }}\right)^{(M-1)}
		\end{aligned},
	\end{equation}
	where $\gamma=\frac{A_s^2}{2(\sigma^2_n + \sigma^2_i)}$ is the SINR of received echo signals.
	The probability of estimated velocity of the target $\underset{j \ne k}{V_e }  = \frac{j \cdot c}{2f_cTM}$ can be derived as
	\begin{equation}\label{equ:Doppler_ind_4}
		\begin{aligned}
			p\left(\underset{j \ne k}{V_e } = \frac{j \cdot c}{2f_cTM}\right) &= e^{-{4\pi^2\gamma k_r^2(m) }} \ll 1
		\end{aligned}.
	\end{equation}
	Proof: See Appendix \ref{app:A}.

	The expectation of estimated velocity $E[V_e]$ can be deduced as
	\begin{equation}\label{equ:Doppler_ind_6}
		\begin{aligned}
			E[V_e] &= \sum_{i=0}^{M-1} V_e \cdot p\left(V_e = \frac{i\cdot c}{2f_cTM} \right) \\
			&\approx V_r \left(1-e^{-{4\pi^2\gamma k_r^2(m) }}\right)^{(M-1)}
		\end{aligned}.
	\end{equation}
	Then, the root mean squared error (RMSE) of velocity can be expressed as
	\begin{equation}\label{equ:Doppler_ind_5}
		\begin{aligned}
			\Gamma_v &= \sqrt{V_r^2 - E[V_e]^2} \\
			&= V_r \sqrt{1-\left(1-e^{-{4\pi^2\gamma k_r^2(m) }}\right)^{2(M-1)}}
		\end{aligned}.
	\end{equation}
	Applying inverse DFT (IDFT) for each column of ${{\bf{S}}_g}$, the range between SBS and the detection target ${R_r}$ can be deduced as follows
	\begin{equation}\label{equ:distance}
		{R_r} \in \left[\frac{{in{d_{s,m}} \cdot c}}{{2{B}}},\frac{{(in{d_{s,m}} + 1) \cdot c}}{{2{B}}}\right),
	\end{equation}
	where $B=N \cdot \Delta f$ is total signal bandwidth, $in{d_{s,m}}$ is the index of the peak of the IDFT of the ${\rm{m}}$-th column of ${{\bf{S}}_g}$. Similarly, RMSE of range can be expressed as
	\begin{equation}\label{equ:Doppler_ind_7}
		\begin{aligned}
			\Gamma_r &= R_r \sqrt{1-\left(1-e^{-{4\pi^2\gamma k_v^2(n) }}\right)^{2(N-1)} }
		\end{aligned}.
	\end{equation}

\begin{figure}[ht]
	\includegraphics[scale=0.49]{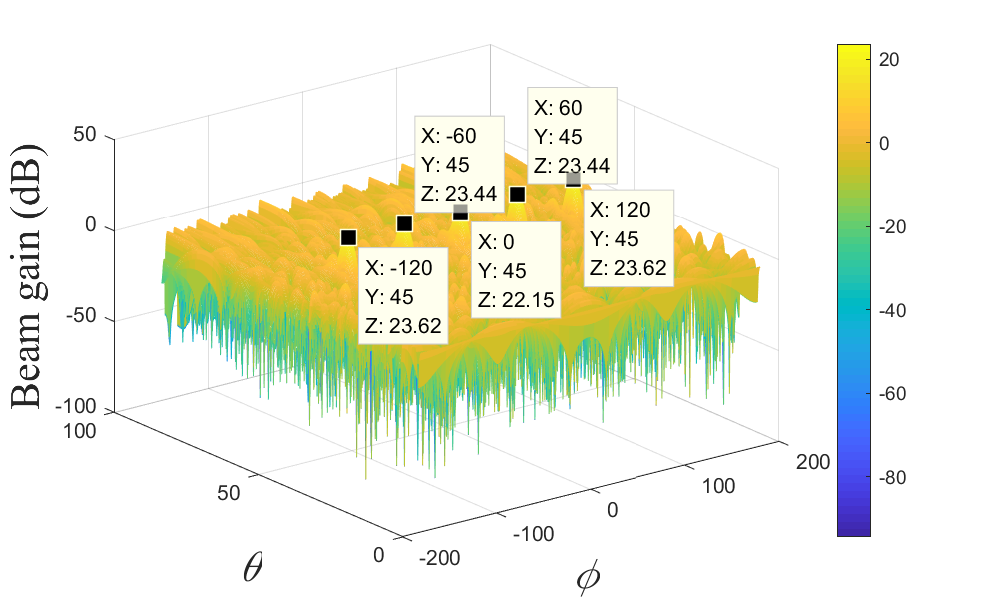}
	\centering
	\caption{Beam gain with direction of ${(-120^o; 45^o)}$, ${(-60^o; 45^o)}$, ${(0^o; 45^o)}$, ${(60^o; 45^o)}$ and ${(120^o; 45^o)}$.}
	\label{fig:SBS_BF_1}
\end{figure}
\begin{figure}[ht]
	\includegraphics[scale=0.49]{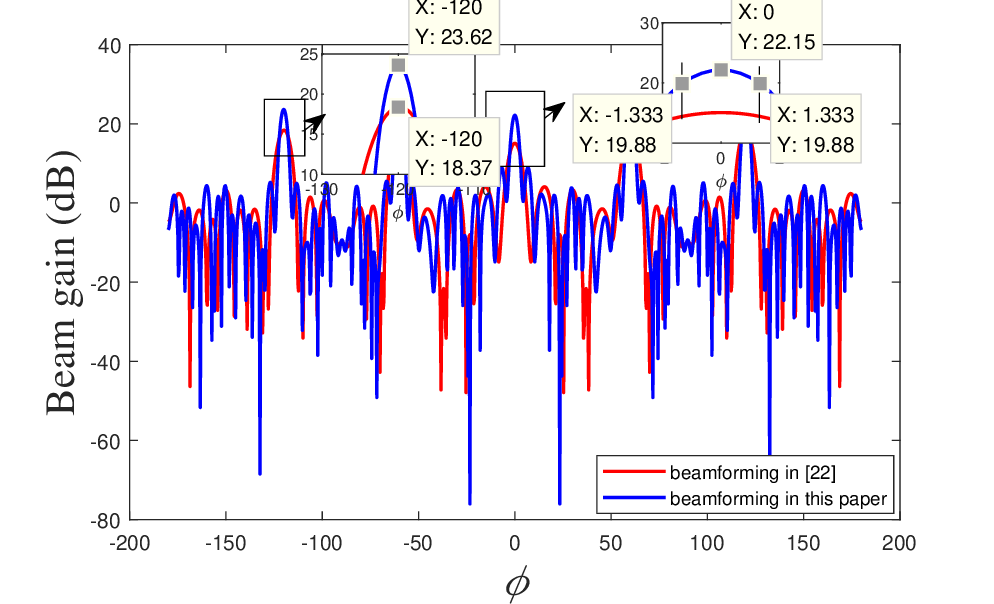}
	\centering
	\caption{Beam gain obtained by two different beamforming algorithms.}
	\label{fig:SBS_BF_2}
\end{figure}
\begin{figure}[ht]
	\includegraphics[scale=0.49]{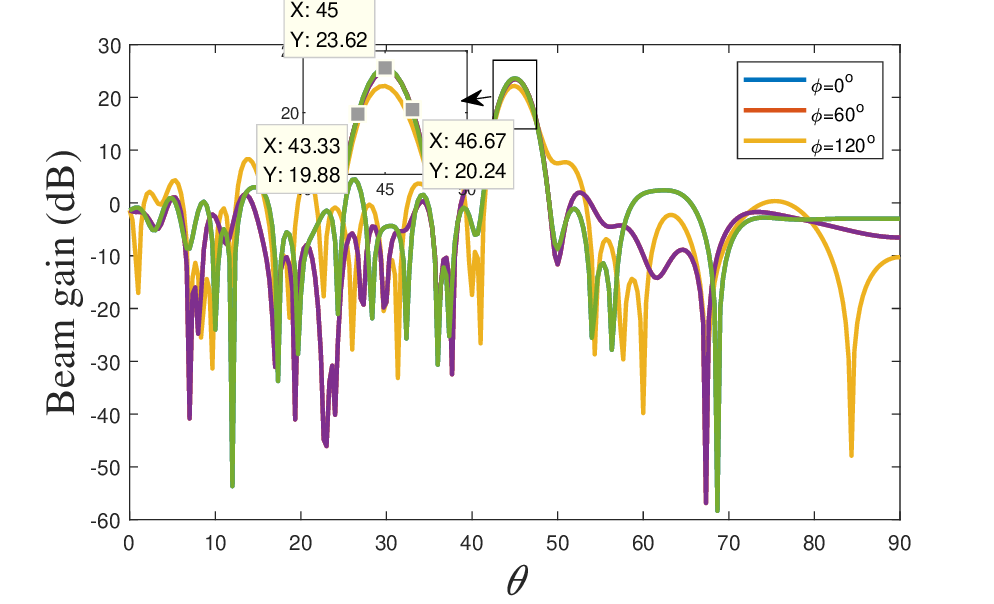}
	\centering
	\caption{Beam gain projection in the direction of $\theta$.}
	\label{fig:SBS_BF_3}
\end{figure}

\subsubsection{Scanning Period of FDB} \label{sec:sensing-Performance-3}
According to \cite{[Radar_dwell]}, the beam dwell duration $\tau $ is 10 ms. The scanning period of FDB, $T_{sc}$ can be derived by the beam dwell duration and the scanning mode of FDB mentioned in section \ref{sec:avoiding-interference}. Details of the derivation to get scanning period are given in algorithm \ref{alg:Framwork_3}.

\section{Performance Simulation of SBS}\label{sec:Numerical results}

In this section, simulations are provided to verify the theoretical analysis in the previous sections.
\subsection{Beamforming of DCB and FDB}

Configurations of FDB, DCB and the antenna array are shown in Table \ref{Parameter:simulation-1}. Beamforming gain with direction of ${\left(-120^o; 45^o\right)}$, ${(-60^o; 45^o)}$, ${(0^o; 45^o)}$, ${(60^o; 45^o)}$ and ${(120^o; 45^o)}$ is shown in Fig. \ref{fig:SBS_BF_1}. Compared with the beamforming algorithm proposed in \cite{[Multibeam]}, the beam gain obtained by the beamforming algorithm proposed in this paper is improved by about 5 dB, just as Fig. \ref{fig:SBS_BF_2} shows. Moreover, Fig. \ref{fig:SBS_BF_2} and Fig. \ref{fig:SBS_BF_3} shows that the width of FDB $(\Delta \theta, \Delta \phi) \approx (2^o,3^o)$.\\
\begin{table}[]
	\caption{Simulation parameters of beamforming.}
	\label{Parameter:simulation-1}
	\begin{tabular}{l|l|l}
		\hline
		\hline
		Items                 & Value               & Meaning of the item     \\ \hline
		${\it{p}}$                       & 16+1                & Antenna layers     \\ \hline
		${2^b}$                          & 16                  & Antennas per layer \\ \hline
		$\beta$ & 0.01 & Regularization factor \\ \hline
		$\left({\phi _{\it{f}}},{\theta _f}\right)$ & $\left({0^o},{45^o}\right)$    & Direction of FDB    \\ \hline
		$\left({\phi _{1}},{\theta _1}\right)$      &  $\left({-120^o},{45^o}\right)$ &                            \\
		$\left({\phi _{2}},{\theta _2}\right)$      & $\left({-60^o},{45^o}\right)$  &  Directions of DCBs \\
		$\left({\phi _{3}},{\theta _3}\right)$      & $\left({60^o},{45^o}\right)$   &                              \\
		$\left({\phi _{4}},{\theta _4}\right)$      & $\left({120^o},{45^o}\right)$  &      \\ \hline
	\end{tabular}
\end{table}

\begin{table}[]
	\caption{Simulation parameters of sensing and communication.}
	\label{Parameter:simulation-2}
	\begin{tabular}{l|l|l}
		\hline
		\hline
		Items                 & Value               & Meaning of the item     \\ \hline
		$\it{h}$ & 10 m & Height of antenna array \\ \hline
		$W_r$ & 20 m & Road width \\ \hline
		$\it{f_c}$ & 24 GHz & Carrier frequency \\ \hline
		${P_{t}} $ & 20 dBm & Total transmit power \\ \hline
		$\rho $ & $0 < \rho < 1$ & Power ratio for sensing \\ \hline
		$\it{k}$ & $1.38*{10^{ - 23}} $ J/K & Boltzmann constant \\ \hline
		${T_{abs}}$ & 290 K & Absolute temperature \\ \hline
		$P_n$ & -94 dBm & Thermal noise \\ \hline
		$P_I$ & -110 dBm & Ground clutter power \\ \hline
		${F_n}$ & 6 dB & \begin{tabular}[c] {@{}l@{}} Noise figure of sensing \\ and communication \end{tabular} \\ \hline
		$\sigma $ & 0.1 ${\rm{m}}^2$ & sensing cross section \\ \hline
		${\it \Gamma}_{\min}$ & 10 dB & \begin{tabular}[c] {@{}l@{}} Minimum required \\ SINR \end{tabular}    \\ \hline
		$\alpha $ & 2.6 & Path loss factor \\ \hline
		$K$ & 10 & Les factor \\ \hline
		$\varepsilon $ & 0.1 & \begin{tabular}[c] {@{}l@{}} Interruption probability \\ threshold \end{tabular}  \\
		\hline
		$\xi_{th} $ & 1,5,10 dB & \begin{tabular}[c] {@{}l@{}} SNR threshold  \\ for communication \end{tabular}  \\
		\hline
		$M$ & 256 & Number of OFDM symbols \\ \hline
		{\color {black} $M_D$} & {\color {black} 2560} & {\color {black} Number of DFT points} \\ \hline
		$N$ & 1024 & Number of subcarriers \\ \hline
		{\color {black} $N_I$} & {\color {black} 10240} & {\color {black} Number of IDFT points} \\ \hline
		${\Delta f}$ & 90.909 kHz & Subcarrier spacing \\ \hline
		${\it{B}}$ & 93.1 MHz & Total signal bandwidth \\ \hline
		$g_{p,s}$ & 54.2 dB & Sensing processing gain \\ \hline
		$g_{p,c}$ & 10 dB & \begin{tabular}[c] {@{}l@{}} Communication \\ processing gain \end{tabular}   \\ \hline
		$g_{r,s}$ & 20 dB & \begin{tabular}[c] {@{}l@{}} Sensing receiving \\ beam gain \end{tabular}  \\ \hline
		$g_{r,c}$ & 6 dB & \begin{tabular}[c] {@{}l@{}} Communication \\ receiving beam gain \end{tabular} \\ \hline
		${{\it{g}}_{t}}$ & 20 dB & Transmitting beam gain \\ \hline
		$R_r$ & 200 m & Range of the target \\ \hline
		$V_r$ & 50 m/s & Velocity of the target \\ \hline
	\end{tabular}
\end{table}

\begin{figure}[ht]
	\includegraphics[scale=0.55]{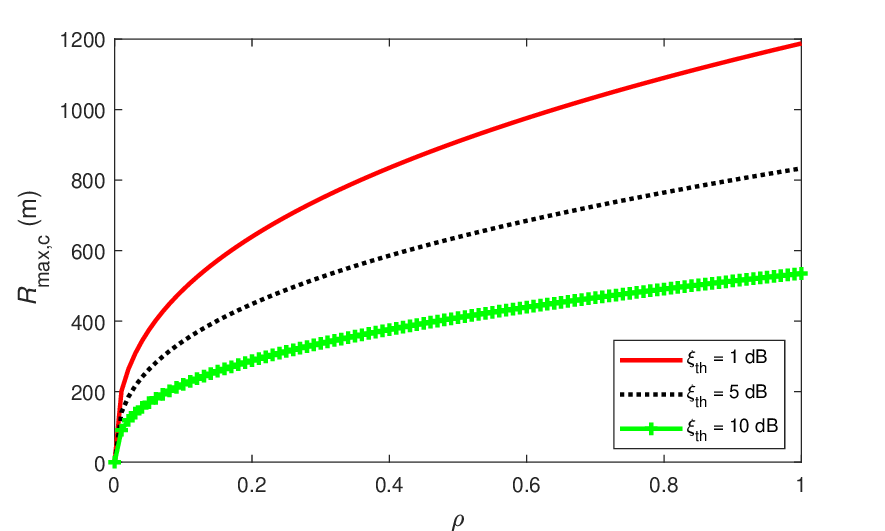}
	\centering
	\caption{Maximum communication range with different $\rho$ in $\xi_{th}=1$ {\rm dB}, $\xi_{th}=5$ {\rm dB} and $\xi_{th}=10$ {\rm dB}.}
	\label{fig:Com_1}
\end{figure}
\begin{figure}[ht]
	\includegraphics[scale=0.55]{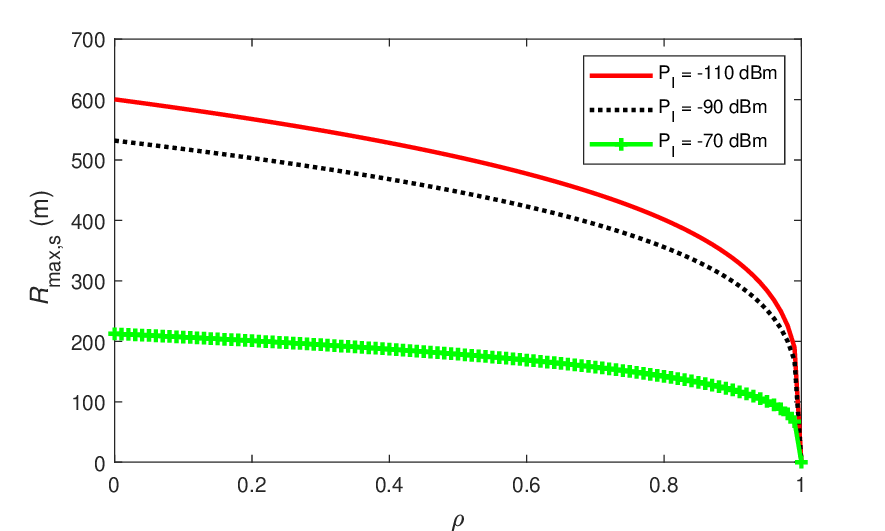}
	\centering
	\caption{Maximum sensing range with different $\rho$ and $P_{I}$.}
	\label{fig:Radar_1}
\end{figure}

\begin{figure}[ht]
	\includegraphics[scale=0.55]{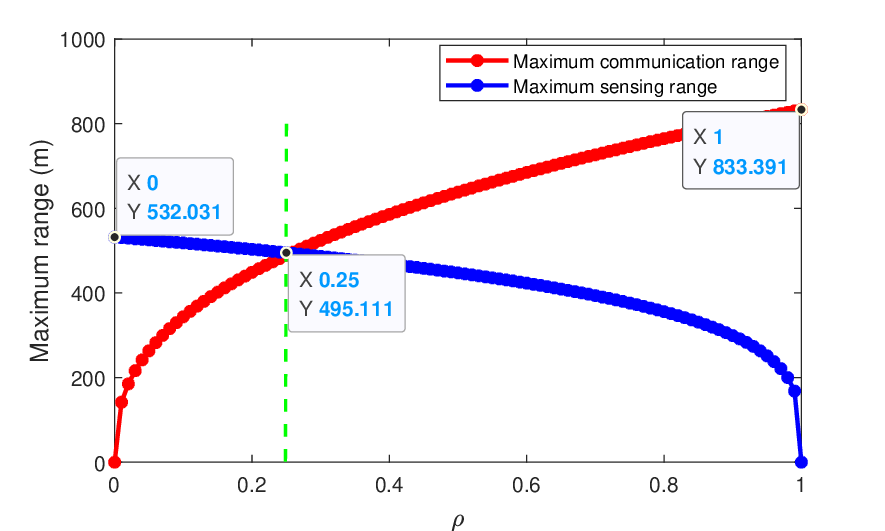}
	\centering
	\caption{{\color{black} Maximum communication and sensing range with different $\rho$ in $\xi_{th}=5$ {\rm dB} and $P_{I}=-90$ {\rm dBm}.}}
	\label{fig:R_sc_2}
\end{figure}

\subsection{Communication and Sensing Performance}
The parameters used in the simulations are shown in Table \ref{Parameter:simulation-2} \cite{[OFDM]} \cite{[parameters_1],[parameters_2],[parameters_3],[24GHz_1],[24GHz_2]}. Due to the high attenuation in radar sensing measurements, and the power limitations of the 24 GHz industrial scientific medical (ISM) regulation of 20 dBm  \cite{[OFDM]}. Therefore, the output power of SBS ${P_{t}}$ = 20 dBm. Based on the current width of the traffic surface and the height of the traffic lights, we set the SBS height and road width to 10 m and 40 m, respectively.

\begin{figure}[ht]
	\includegraphics[scale=0.55]{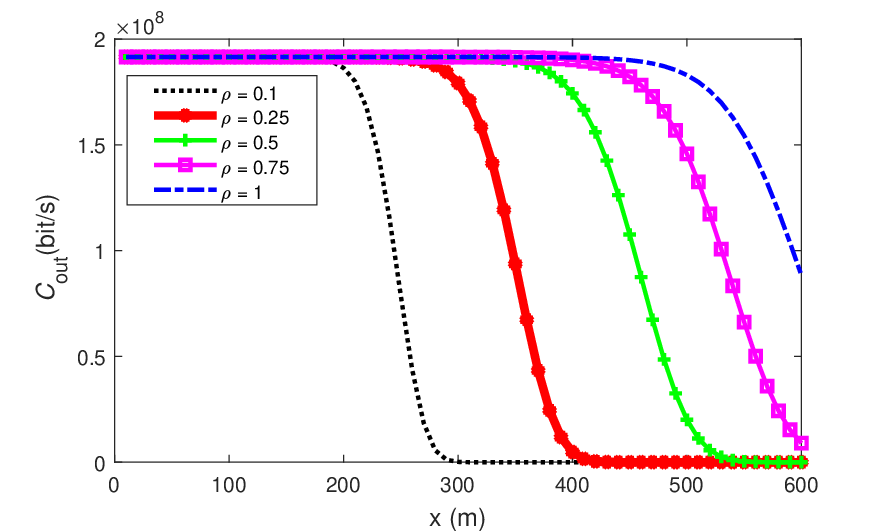}
	\centering
	\caption{The outage capacity with different communication distance in $\xi_{th}=5$ {\rm dB}.}
	\label{fig:SBS_Tc_2}
\end{figure}

\subsubsection{Maximum Communication and Sensing Range}
According to \eqref{equ:Radar_range} and \eqref{equ:Communication_range}, both the maximum sensing range $R_{\max,s}$ and the maximum communication range $R_{\max,c}$ depend on the power ratio for communication $\rho$. With the increase of $\rho$, the power used for communication will be larger, $R_{\max,c}$ will be larger, just as Fig. \ref{fig:Com_1} shows. Conversely, $R_{\max,s}$ will be smaller just as Fig. \ref{fig:Radar_1} shows. Moreover, the SNR threshold for communication $\xi_{th}$ is another important factor affecting $R_{\max,c}$. As Fig. \ref{fig:Com_1} shows, $R_{\max,c}$ decreases as $\xi_{th}$ increases. As for $R_{\max,s}$, the ground clutter power $P_I$ is an important influence factor. As Fig. \ref{fig:Radar_1} shows, the higher $P_I$ is, the smaller $R_{\max,s}$ is.

To analyze the tradeoff between the communication and sensing performance in terms of range, Fig. \ref{fig:R_sc_2} show the maximum communication and sensing range with different $\rho$ in $\xi_{th}= 1$ dB and $P_{I} = -100$ dBm. When $\rho=0$, all power of SBS is used for sensing, and $R_{\max,s}$ can reach about 532 m. When $\rho=1$, all power of SBS is used for communication, and $R_{\max,c}$ can reach about 833 m. Moreover, we can adjust $\rho$ to 0.25 so that communication and sensing get the same range of 495.111 m.

\begin{figure}[ht]
	\includegraphics[scale=0.55]{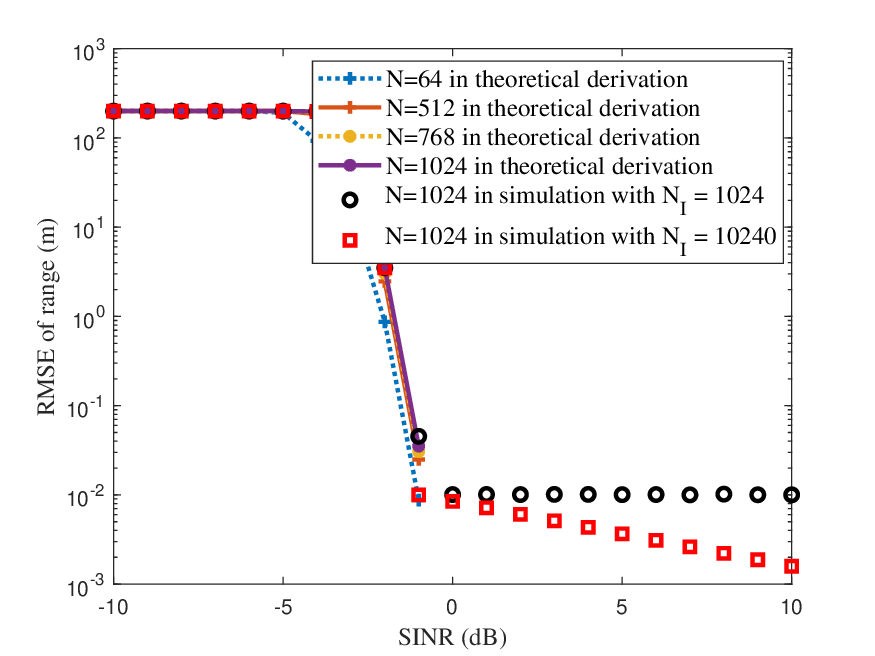}
	\centering
	\caption{\color {black} Accuracy of the range of target with different SINR.}
	\label{fig:SBS_ACC_R}
\end{figure}
\begin{figure}[ht]
	\includegraphics[scale=0.55]{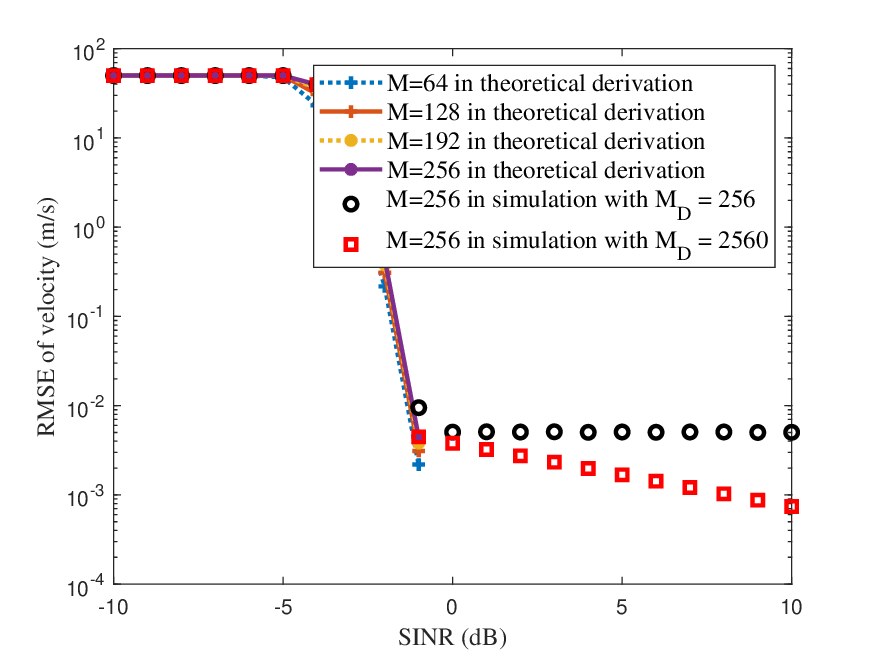}
	\centering
	\caption{\color {black} Accuracy of the velocity of target with different SINR.}
	\label{fig:SBS_ACC_V}
\end{figure}

\subsubsection{Communication Outage Capacity}
Based on \eqref{equ:OutThrouResult}, the outage capacity is influenced by the communication distance and beamwidth. In theory, the path loss grows with the communication distance. Severe path loss further decreases the outage capacity. As Fig. \ref{fig:SBS_Tc_2} shows, the outage capacity decreases with the communication distance increases.

\subsubsection{RMSE of Range and Velocity}
	As Fig. \ref{fig:SBS_ACC_R} and Fig. \ref{fig:SBS_ACC_V} shows, the simulation accuracy of the range and velocity is close to the theoretical derivation, which verifies the derivation. {\color{black} Each simulation in this section is calculated over 5000 Monte Carlo trials.}
	Moreover, the range accuracy performance of the target increases with the increase of the number of subcarriers $N$, the velocity accuracy performance of the target increases with the increase of the number of symbols $M$.

\subsubsection{Scanning Period}

\begin{figure}[ht]
	\includegraphics[scale=0.55]{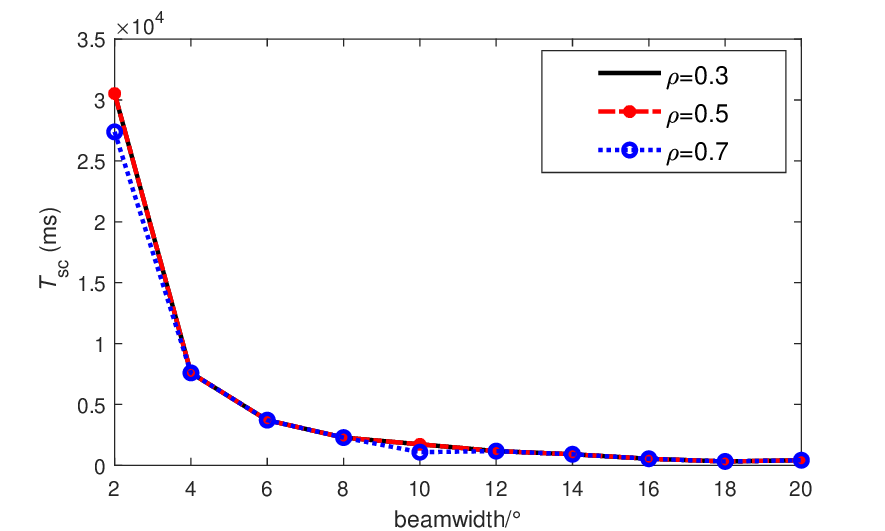}
	\centering
	\caption{Scanning period with different beamwidth and $\rho$.}
	\label{fig:SBS_scan_period_2}
\end{figure}

The scanning period can be derived with the help of the scanning mode mentioned in section \ref{sec:avoiding-interference}. In theory, as $\rho$ increases, the smaller the power used for sensing, the smaller the maximum sensing range $R_{\max,s}$, and the smaller the scanning period should be. However, as $R_{\max,s}$ does not change much, the scanning period does not change much, just as Fig. \ref{fig:SBS_scan_period_2} shows. Moreover, with the increase of beamwidth, the range of beam coverage is larger, and the scanning period is also reduced.

\section{Conclusion}\label{sec:Conclusion}
In this paper, we study the SBS system that has great potential to improve the safety for vehicles and pedestrians on roads. Key technologies of SBS are studied, including beamforming algorithm, beam scanning scheme, and interference cancellation algorithm.
The interference of single SBS system contains multiuser interference and duplex interference between sensing and communication.
We design a new signal frame structure to cancel the duplex interference, and propose the TSF-DMA algorithm to cancel the multiuser interference.
The communication and sensing performance of the SBS system is simulated and analyzed. According to the simulation results, the allocation of power resources can affect the performance of sensing and communication, including the communication and sensing range. We can adjust the communication and sensing range by the power ratio for communication $\rho$.

\appendix
	\section{Probability of $V_e = \frac{k \cdot c}{2f_cTM}$} \label{app:A}
	Under the condition of $V_e = \frac{k \cdot c}{2f_cTM}$, we can get
	\begin{equation}\label{equ:APP_A1}
		k_v(m)=e^{j2\pi\frac{mk}{M}}.
	\end{equation}
	Then $v(k,n)$ can be deduced as
	\begin{equation}\label{equ:APP_A2}
		\begin{aligned}
			v(k,n) &= {A_s}{k_r(n)} \sum_{m=0}^{M-1} e^{j2\pi\frac{mk}{M}}e^{-j2\pi\frac{mk}{M}} + \\ & \sum_{m=0}^{M-1}N_Ie^{-j2\pi\frac{mk}{M}} \\
			&= {A_s}{k_r(n)}M + \sum_{m=0}^{M-1}N_Ie^{-j2\pi\frac{mk}{M}}.
		\end{aligned}
	\end{equation}
	Similarly, $v(i,n)$ can be deduced as
	\begin{equation}\label{equ:APP_A3}
		\begin{aligned}
			v(i,n) &= {A_s}{k_r(n)} \begin{Vmatrix} \sum_{m=0}^{M-1} e^{j2\pi\frac{mk}{M}}e^{-j2\pi\frac{im}{M}}\end{Vmatrix} + \\ & \sum_{m=0}^{M-1}N_Ie^{-j2\pi\frac{im}{M}} \\
			&= {A_s}{k_r(n)} \begin{Vmatrix} \frac{1-e^{j2\pi(k-i)}}{1-e^{j2\pi\frac{k-i}{M}}}\end{Vmatrix} + \sum_{m=0}^{M-1}N_Ie^{-j2\pi\frac{im}{M}} \\
			&=
			{A_s}{k_r(n)}\frac{\sin( (k-i)\pi)}{\sin( \frac{k-i}{M}\pi)} + \sum_{m=0}^{M-1}N_Ie^{-j2\pi\frac{im}{M}} \\
			&\approx
			{A_s}{k_r(n)} \frac{\sin( (k-i)\pi)}{ \frac{k-i}{M}\pi} + \sum_{m=0}^{M-1}N_Ie^{-j2\pi\frac{im}{M}} \\
			&=
			{A_s}{k_r(n)}M {\rm sinc}(k-i) + \sum_{m=0}^{M-1}N_Ie^{-j2\pi\frac{im}{M}},
		\end{aligned}
	\end{equation}
	where ${\rm sinc}(x) = \frac{\sin(\pi x)}{\pi x}$.
	Then $p\left(v(k,n) \ge v(i,n)\right)$ can be simplified as
	\begin{equation}\label{equ:APP_A4}
		\begin{aligned}
			p\left(C \ge \max \left\{\sum_{m=0}^{M-1}N_I \left(e^{-j2\pi\frac{im}{M}}-e^{-j2\pi\frac{mk}{M}}\right) \right\} \right),
		\end{aligned}
	\end{equation}
	where
	\begin{equation}\label{equ:APP_A5}
		C={A_s}{k_r(n)}M\left(1-{\rm sinc}(k-i)\right).
	\end{equation}

	Since $N_I = A_n + A_I = A_{ni}e^{j\varphi_{ni}}$, the real and image of $N_I$ follow gaussian distribution ${{\mathcal{CN}}\rm{(0,}}\sigma _n^2+\sigma _i^2{\rm{)}}$. The phase $\varphi_{ni}$ follows uniform distribution. The amplitude $A_{ni}$ follows Rayleigh distribution, the probability density function can be expressed as
	\begin{equation}\label{equ:rayleigh}
		{\it{f}_{NI}(A_{ni})} =\left\{
		\begin{aligned}
			\frac{A_{ni}}{\sigma^2_n + \sigma^2_i}e^{-\frac{A_{ni}^2}{2 \left(\sigma^2_n + \sigma^2_i\right)}}, & \quad A_{ni} \ge 0 \\
			0, & \quad A_{ni} \le 0
		\end{aligned}.
		\right.
	\end{equation}
	The right side of \eqref{equ:APP_A4} can be expressed as
	\begin{equation}\label{equ:APP_A7}
		\begin{aligned}
			\sum_{m=0}^{M-1} A_{ni}e^{j\varphi_{ni}} \left(e^{-j2\pi\frac{im}{M}}-e^{-j2\pi\frac{mk}{M}}\right) , & \quad x \ge 0
		\end{aligned}.
	\end{equation}
	Similarly, we just need to consider the case of  $\varphi_{ni} = 2\pi\frac{im}{M}$. Then the maximum of the right side of \eqref{equ:APP_A4} is approximately equal to $\frac{M}{2\pi}A_{ni} $.

	Then $p\left(v(k,n) \ge v(i,n)\right)$ can be expressed as
	\begin{equation}\label{equ:APP_A9}
		\begin{aligned}
			&= p\left(2\pi A_sk_r(n)(1-{\rm sinc}(k-i)) \ge A_{ni}\right) \\
			&= \int_{0}^{2\pi A_s k_r(n) (1-{\rm sinc}(k-i))} {\frac{x}{\sigma^2_n + \sigma^2_i}e^{-\frac{x^2}{2(\sigma^2_n + \sigma^2_i)}}} dx \\
			&= 1-e^{-\frac{2\pi^2 A_s^2 k_r^2(m) \left(1-{\rm sinc}(k-i)\right)^2}{\sigma^2_n + \sigma^2_i}} \\
			&\approx
			1-e^{-\frac{2\pi^2 A_s^2 k_r^2(m)}{\sigma^2_n + \sigma^2_i}}
		\end{aligned}.
	\end{equation}
	Then
	\begin{equation}\label{equ:APP_A10}
		\begin{aligned}
			p\left(V_e = \frac{k \cdot c}{2f_cTM}\right) &= \prod_{i=0,i\ne k}^{M-1} p\left(v(k,n) \ge v(i,n)\right) \\
			&= \left(1-e^{-\frac{2\pi^2 A_s^2 k_r^2(m)}{\sigma^2_n + \sigma^2_i}}\right)^{(M-1)}
		\end{aligned}.
	\end{equation}

	The probability of estimated velocity of the target $\underset{j \ne k}{V_e }  = \frac{j \cdot c}{2f_cTM}$ can be expressed as
	\begin{equation}\label{equ:Doppler_ind_B1}
		\begin{aligned}
			p\left(\underset{j \ne k}{V_e } = \frac{j \cdot c}{2f_cTM}\right)
			&= \prod_{i=0,i\ne j}^{M-1} p\left(v(j,n) >  v(i,n)\right) \\
			&\approx p(v(j,n) \ge v(k,n)) \\
			&= 1 - p(v(j,n) \le v(k,n)) \\
			&= e^{-\frac{2\pi^2 A_s^2 k_r^2(m)}{\sigma^2_n + \sigma^2_i}}
		\end{aligned}.
	\end{equation}

\bibliographystyle{IEEEtran}
\bibliography{reference}

\begin{thebibliography}{10}
\providecommand{\url}[1]{#1}
\csname url@samestyle\endcsname
\providecommand{\newblock}{\relax}
\providecommand{\bibinfo}[2]{#2}
\providecommand{\BIBentrySTDinterwordspacing}{\spaceskip=0pt\relax}
\providecommand{\BIBentryALTinterwordstretchfactor}{4}
\providecommand{\BIBentryALTinterwordspacing}{\spaceskip=\fontdimen2\font plus
\BIBentryALTinterwordstretchfactor\fontdimen3\font minus
  \fontdimen4\font\relax}
\providecommand{\BIBforeignlanguage}[2]{{%
\expandafter\ifx\csname l@#1\endcsname\relax
\typeout{** WARNING: IEEEtran.bst: No hyphenation pattern has been}%
\typeout{** loaded for the language `#1'. Using the pattern for}%
\typeout{** the default language instead.}%
\else
\language=\csname l@#1\endcsname
\fi
#2}}
\providecommand{\BIBdecl}{\relax}
\BIBdecl

\bibitem{[Driverless_Future]}
J.~{Rowley}, A.~{Liu}, S.~{Sandry}, J.~{Gross}, M.~{Salvador}, C.~{Anton}, and
  C.~{Fleming}, ``Examining the driverless future: An analysis of human-caused
  vehicle accidents and development of an autonomous vehicle communication
  testbed,'' in \emph{2018 Systems and Information Engineering Design Symposium
  (SIEDS)}, 2018, pp. 58--63.

\bibitem{[C-ITS]}
L.~{Chen} and C.~{Englund}, ``Cooperative intersection management: A survey,''
  \emph{IEEE Transactions on Intelligent Transportation Systems}, vol.~17,
  no.~2, pp. 570--586, 2016.

\bibitem{[RSU]}
W.~{Liu}, S.~{Muramatsu}, and Y.~{Okubo}, ``Cooperation of {V2I/P2I}
  communication and roadside radar perception for the safety of vulnerable road
  users,'' in \emph{2018 16th International Conference on Intelligent
  Transportation Systems Telecommunications (ITST)}, 2018, pp. 1--7.

\bibitem{[VMSE]}
J.~{Wang}, J.~{Liu}, and N.~{Kato}, ``Networking and communications in
  autonomous driving: A survey,'' \emph{IEEE Communications Surveys \&
  Tutorials}, vol.~21, no.~2, pp. 1243--1274, 2019.

\bibitem{[RSU_VMSE]}
D.~{Kim}, Y.~{Velasco}, W.~{Wang}, R.~N. {Uma}, R.~{Hussain}, and S.~{Lee}, ``A
  new comprehensive {RSU} installation strategy for cost-efficient {VANET}
  deployment,'' \emph{IEEE Transactions on Vehicular Technology}, vol.~66,
  no.~5, pp. 4200--4211, 2017.

\bibitem{[JWJ]}
W.~Jiang, Z.~Wei, B.~Li, Z.~Feng, and Z.~Fang, ``Improve radar sensing
  performance of multiple roadside units cooperation via space registration,''
  \emph{IEEE Transactions on Vehicular Technology}, vol.~71, no.~10, pp.
  10\,975--10\,990, Oct. 2022.

\bibitem{[RSU_Energy]}
A.~{Khezrian}, T.~D. {Todd}, G.~{Karakostas}, and M.~{Azimifar},
  ``Energy-efficient scheduling in green vehicular infrastructure with multiple
  roadside units,'' \emph{IEEE Transactions on Vehicular Technology}, vol.~64,
  no.~5, pp. 1942--1957, May. 2015.

\bibitem{[JRC_1]}
Z.~{Feng}, Z.~{Fang}, Z.~{Wei}, X.~{Chen}, Z.~{Quan}, and D.~{Ji}, ``Joint
  radar and communication: A survey,'' \emph{China Communications}, vol.~17,
  no.~1, pp. 1--27, Jan. 2020.

\bibitem{[spectrum_limit]}
Z.~{Chen}, X.~{Ma}, B.~{Zhang}, Y.~{Zhang}, Z.~{Niu}, N.~{Kuang}, W.~{Chen},
  L.~{Li}, and S.~{Li}, ``A survey on terahertz communications,'' \emph{China
  Communications}, vol.~16, no.~2, pp. 1--35, Feb. 2019.

\bibitem{[JRC_AV]}
D.~{Ma}, N.~{Shlezinger}, T.~{Huang}, Y.~{Liu}, and Y.~C. {Eldar}, ``Joint
  radar-communication strategies for autonomous vehicles: Combining two key
  automotive technologies,'' \emph{IEEE Signal Processing Magazine}, vol.~37,
  no.~4, pp. 85--97, July. 2020.

\bibitem{[6G]}
W.~{Saad}, M.~{Bennis}, and M.~{Chen}, ``A vision of {6G} wireless systems:
  Applications, trends, technologies, and open research problems,'' \emph{IEEE
  Network}, vol.~34, no.~3, pp. 134--142, May/June. 2020.

\bibitem{[JRC_waveform]}
S.~H. {Dokhanchi}, M.~R. {Bhavani Shankar}, M.~{Alaee-Kerahroodi},
  T.~{Stifter}, and B.~{Ottersten}, ``Adaptive waveform design for automotive
  joint radar-communications system,'' in \emph{ICASSP 2019 - 2019 IEEE
  International Conference on Acoustics, Speech and Signal Processing
  (ICASSP)}, 2019, pp. 4280--4284.

\bibitem{[OFDM]}
C.~{Sturm} and W.~{Wiesbeck}, ``Waveform design and signal processing aspects
  for fusion of wireless communications and radar sensing,'' \emph{Proceedings
  of the IEEE}, vol.~99, no.~7, pp. 1236--1259, July. 2011.

\bibitem{[LFM]}
Y.~{Zhang}, Q.~{Li}, L.~{Huang}, K.~{Dai}, and J.~{Song}, ``Waveform design for
  joint radar-communication with nonideal power amplifier and outband
  interference,'' in \emph{2017 IEEE Wireless Communications and Networking
  Conference (WCNC), San Francisco, CA}, 2017, pp. 1--6.

\bibitem{[PMCW]}
S.~H. {Dokhanchi}, B.~S. {Mysore}, K.~V. {Mishra}, and B.~{Ottersten}, ``A
  mmwave automotive joint radar-communications system,'' \emph{IEEE
  Transactions on Aerospace and Electronic Systems}, vol.~55, no.~3, pp.
  1241--1260, June. 2019.

\bibitem{[OFDM_2]}
C.~{Sturm}, T.~{Zwick}, W.~{Wiesbeck}, and M.~{Braun}, ``Performance
  verification of symbol-based {OFDM} radar processing,'' in \emph{2010 IEEE
  Radar Conference}, Washington. DC, 2010, pp. 60--63.

\bibitem{[OFDM_advantage]}
Y.~{Liu}, G.~{Liao}, J.~{Xu}, Z.~{Yang}, and Y.~{Zhang}, ``Adaptive {OFDM}
  integrated radar and communications waveform design based on information
  theory,'' \emph{IEEE Communications Letters}, vol.~21, no.~10, pp.
  2174--2177, Oct. 2017.

\bibitem{[OFDM_3]}
X.~{Tian}, T.~{Zhang}, Q.~{Zhang}, and Z.~{Song}, ``Waveform design and
  processing in {OFDM} based radar-communication integrated systems,'' in
  \emph{2017 IEEE/CIC International Conference on Communications in China
  (ICCC)}, Qingdao, 2017, pp. 1--6.

\bibitem{[IIE-1]}
Y.~L. {Sit}, C.~{Sturm}, and T.~{Zwick}, ``One-stage selective interference
  cancellation for the {OFDM} joint radar-communication system,'' in \emph{2012
  The 7th German Microwave Conference, Ilmenau}, 2012, pp. 1--4.

\bibitem{[IIE-2]}
Y.~L. {Sit}, B.~{Nuss}, S.~{Basak}, M.~{Orzol}, and T.~{Zwick}, ``Demonstration
  of interference cancellation in a multiple-user access {OFDM MIMO}
  radar-communication network using {USRPs},'' in \emph{2016 IEEE MTT-S
  International Conference on Microwaves for Intelligent Mobility (ICMIM), San
  Diego, CA}, 2016, pp. 1--4.

\bibitem{[IIE-3]}
Y.~L. {Sit}, B.~{Nuss}, and T.~{Zwick}, ``On mutual interference cancellation
  in a {MIMO OFDM} multiuser radar-communication network,'' \emph{IEEE
  Transactions on Vehicular Technology}, vol.~67, no.~4, pp. 3339--3348, April.
  2018.

\bibitem{[Multibeam]}
J.~A. {Zhang}, X.~{Huang}, Y.~J. {Guo}, J.~{Yuan}, and R.~W. {Heath},
  ``Multibeam for joint communication and radar sensing using steerable analog
  antenna arrays,'' \emph{IEEE Transactions on Vehicular Technology}, vol.~68,
  no.~1, pp. 671--685, Jan. 2019.

\bibitem{[rice_vehicular]}
C.~A. {Gutierrez} and C.~A. {Gómez-Vega}, ``A {Non-WSSUS} rice fading channel
  model for vehicular communications,'' in \emph{2018 IEEE 10th Latin-American
  Conference on Communications (LATINCOM)}, 2018, pp. 1--6.

\bibitem{[RICE]}
M.~K. {Simon} and M.~{Alouini}, ``A unified approach to the performance
  analysis of digital communication over generalized fading channels,''
  \emph{Proceedings of the IEEE}, vol.~86, no.~9, pp. 1860--1877, Sep. 1998.

\bibitem{[Xinyuan]}
X.~Yuan, Z.~Feng, W.~Xu, W.~Ni, J.~A. Zhang, Z.~Wei, and R.~P. Liu, ``Capacity
  analysis of {UAV} communications: Cases of random trajectories,'' \emph{IEEE
  Transactions on Vehicular Technology}, vol.~67, no.~8, pp. 7564--7576, Aug.
  2018.

\bibitem{[Xinyuan-4]}
M.~M. Azari, F.~Rosas, K.-C. Chen, and S.~Pollin, ``Ultra reliable {UAV}
  communication using altitude and cooperation diversity,'' \emph{IEEE
  Transactions on Communications}, vol.~66, no.~1, pp. 330--344, Jan. 2018.

\bibitem{[Xinyuan-4-31]}
M.~K. Simon and M.-S. Alouini, \emph{Digital Communication Over Fading
  Channels}, vol. 95. Hoboken, NJ, USA: Wiley, 2005.

\bibitem{[ant_design]}
R.~{Mailloux}, ``Phased array antenna handbook, third edition,'' \emph{Phased
  Array Antenna Handbook, Third Edition, Artech}, 2017.

\bibitem{[UCA]}
P.~Ioannides and C.~Balanis, ``Uniform circular and rectangular arrays for
  adaptive beamforming applications,'' \emph{IEEE Antennas and Wireless
  Propagation Letters}, vol.~4, pp. 351--354, Sept. 2005.

\bibitem{[duplex_interference]}
A.~Liu, W.~Sheng, and T.~Riihonen, ``Per-antenna self-interference cancellation
  beamforming design for digital phased array,'' \emph{IEEE Signal Processing
  Letters}, vol.~29, pp. 2442--2446, Nov. 2022.

\bibitem{[OFDMA]}
S.~{Yang}, \emph{{OFDMA} System Analysis and Design}, Artech. 2010.

\bibitem{[Outage_Pro]}
S.~{Zhalehpour}, M.~{Uysal}, O.~A. {Dobre}, and T.~{Ngatched}, ``Outage
  capacity and throughput analysis of multiuser {FSO} systems,'' in \emph{2015
  IEEE 14th Canadian Workshop on Information Theory (CWIT)}, 2015, pp.
  143--146.

\bibitem{azari2018ultra}
M.~K. {Simon} and M.-S. {Alouini}, \emph{Outage Performance of Multiuser
  Communication Systems}, 2005, pp. 638--680.

\bibitem{[Q_func]}
D.~A. {Shnidman}, ``The calculation of the probability of detection and the
  generalized marcum {Q}-function,'' \emph{IEEE Transactions on Information
  Theory}, vol.~35, no.~2, pp. 389--400, March. 1989.

\bibitem{[OFDM_resolution]}
J.~{Li}, S.~{An}, J.~{An}, H.~{Zirath}, and Z.~S. {He}, ``{OFDM} radar range
  accuracy enhancement using fractional fourier transformation and phase
  analysis techniques,'' \emph{IEEE Sensors Journal}, vol.~20, no.~2, pp.
  1011--1018, 2020.

\bibitem{[rayleigh]}
X.~{Wang}, H.~{Wang}, S.~{Yan}, L.~{Li}, and C.~{Meng}, ``Simulation for
  surveillance radar ground clutter at low grazing angle,'' in \emph{2012
  International Conference on Image Analysis and Signal Processing}, Hangzhou,
  2012, pp. 1--4.

\bibitem{[compound-gaussian]}
E.~{Ollila}, D.~E. {Tyler}, V.~{Koivunen}, and H.~V. {Poor},
  ``Compound-gaussian clutter modeling with an inverse gaussian texture
  distribution,'' \emph{IEEE Signal Processing Letters}, vol.~19, no.~12, pp.
  876--879, 2012.

\bibitem{[Radar_dwell]}
M.~{Scharrenbroich} and M.~{Zatman}, ``Joint radar-communications resource
  management,'' in \emph{2016 IEEE Radar Conference (RadarConf)}, Philadelphia,
  PA, 2016, pp. 1--6.

\bibitem{[parameters_1]}
E.~{Markin}, \emph{Principles of modern radar missile seekers}, 2022.

\bibitem{[parameters_2]}
H.~{Boeglen}, B.~{Hilt}, P.~{Lorenz}, J.~{Ledy}, A.~{Poussard}, and
  R.~{Vauzelle}, ``A survey of {V2V} channel modeling for {VANET}
  simulations,'' in \emph{2011 Eighth International Conference on Wireless
  On-Demand Network Systems and Services, Bardonecchia}, 2011, pp. 117--123.

\bibitem{[parameters_3]}
V.~{Lakshmanan}, ``Overview of radar data compression,'' \emph{Satellite Data
  Compression, Communications, and Archiving {III}}, vol. 6683, Sep. 2007.

\bibitem{[24GHz_1]}
C.~{Sturm}, T.~{Zwick}, and W.~{Wiesbeck}, ``An {OFDM} system concept for joint
  radar and communications operations,'' in \emph{VTC Spring 2009 - IEEE 69th
  Vehicular Technology Conference}, Barcelona, 2009, pp. 1--5.

\bibitem{[24GHz_2]}
M.~{Braun}, C.~{Sturm}, A.~{Niethammer}, and F.~K. {Jondral}, ``Parametrization
  of joint {OFDM}-based radar and communication systems for vehicular
  applications,'' in \emph{2009 IEEE 20th International Symposium on Personal,
  Indoor and Mobile Radio Communications}, Tokyo, 2009.

\end{thebibliography}
\newpage

\end{document}